\documentclass[11pt]{article}
\usepackage[T1]{fontenc}
\usepackage{lmodern,microtype,amsmath,amssymb,amsthm,mathtools,mathrsfs}
\usepackage[letterpaper,margin=.88in]{geometry}
\newcommand{\cV}{\mathcal V}

\newcommand{\Fix}{\operatorname{Fix}}
\usepackage{enumitem,needspace,longtable,array}
\usepackage[hidelinks]{hyperref}
\usepackage{bookmark}
\makeatletter
\newcommand{\resultbookmark}[1]{\bookmark[level=3,dest=\@currentHref]{\@currentlabel\space #1}}
\makeatother
\newcommand{\proofplace}[1]{\hyperref[#1]{\ref*{#1} (p.~\pageref*{#1})}}
\hypersetup{pdftitle={Real-analytic actions of higher-rank lattices preserving an ergodic volume},pdfauthor={Federico Rodriguez Hertz; Zhiren Wang}}
\numberwithin{equation}{section}
\newtheorem{theorem}[equation]{Theorem}
\newtheorem{lemma}[equation]{Lemma}
\newtheorem{proposition}[equation]{Proposition}

\theoremstyle{definition}\newtheorem{definition}[equation]{Definition}

\theoremstyle{remark}
\counterwithin*{section}{part}
\renewcommand*{\thesection}{\ifnum\value{part}>0 \thepart.\fi\arabic{section}}
\renewcommand*{\thesubsection}{\thesection.\arabic{subsection}}
\renewcommand*{\theequation}{\ifnum\value{subsection}>0 \thesubsection\else\thesection\fi.\arabic{equation}}

\newcommand{\R}{\mathbb R}
\newcommand{\id}{\operatorname{id}}\newcommand{\GL}{\operatorname{GL}}
\newcommand{\Hom}{\operatorname{Hom}}\newcommand{\Sym}{\operatorname{Sym}}
\newcommand{\Lie}{\operatorname{Lie}}\newcommand{\Ad}{\operatorname{Ad}}

\newcommand{\Gr}{\operatorname{Gr}}

\newcommand{\Z}{\mathbb Z}

\title{Real-analytic actions of higher-rank lattices\\preserving an ergodic volume}
\author{Federico Rodriguez Hertz \and Zhiren Wang}
\date{}
\makeatletter
\renewcommand{\@maketitle}{%
  \begin{center}
    {\LARGE \@title\par}\vskip 1em
    {\large\begin{tabular}[t]{c}\@author\end{tabular}\par}%
  \end{center}\vskip .6em}
\makeatother

\newcommand{\vol}{\operatorname{vol}}

\newcommand{\manuscriptpart}[1]{%
 \refstepcounter{part}\setcounter{section}{0}\setcounter{subsection}{0}%
 \section*{Part \thepart.\ #1}%
 \addcontentsline{toc}{part}{Part \thepart.\ #1}}
\makeatletter
\renewcommand{\l@section}{\@dottedtocline{1}{0em}{3.8em}}
\renewcommand{\l@subsection}{\@dottedtocline{2}{1.5em}{4.6em}}
\makeatother
\newtheorem{claim}[equation]{Claim}
\newenvironment{claimproof}{\begin{proof}[Proof of the claim]}{\end{proof}}

\begin{document}

\maketitle
\begin{abstract}
We classify real-analytic actions of higher-rank lattices that preserve an ergodic smooth volume on a compact manifold. On an invariant open dense set of full volume, the action is analytically conjugate to the restriction of an affine action on finitely many copies of one double quotient. The conjugacy carries Haar probability to the given volume. Analytic germ superrigidity supplies local coordinates with measurable dependence on the center. The conditional measures on coarse Lyapunov manifolds determine the law of a displaced center. A two-point joining pairs charts whose evaluations agree there and gives the change of coordinates. The resulting analytic frame bundle carries a coframe whose vector fields integrate to the homogeneous model. This gives the open homogeneous classification in Zimmer's program.
\end{abstract}
\section{The homogeneous model}\label{sec:homogeneous}
Let $G$ be a connected semisimple real linear Lie group with
$\operatorname{rank}_{\R}G_i\ge2$ for every noncompact simple factor $G_i$.
Let a lattice $\Gamma<G$ act by real-analytic diffeomorphisms on a compact
real-analytic manifold $M$ without boundary. Let $m$ be a $\Gamma$-invariant
ergodic probability given in local coordinates by $dm=a(x)\,dx$, $a\in C^\infty$, $a>0$.

The smooth locus of a real-analytic quotient consists of the points having a neighborhood analytically isomorphic, with its quotient analytic structure, to an open subset of a Euclidean space. For a Lie group $L$, a lattice $\Delta<L$ and a compact subgroup $H<L$, write $m_C$ for the image of Haar probability on $\Delta\backslash L$ under $\Delta g\mapsto\Delta gH$ in $C=\Delta\backslash L/H$. An affine automorphism of $C$ has the form $\Delta gH\mapsto\Delta b\alpha(g)H$, with $\alpha\in\operatorname{Aut}(L)$, $\alpha(H)=H$ and $b\alpha(\Delta)b^{-1}=\Delta$.
Given a finite-index subgroup $\Gamma_C<\Gamma$ acting affinely on $C$, put $s=[\Gamma:\Gamma_C]$ and choose representatives $t_1=e,\ldots,t_s$ of $\Gamma/\Gamma_C$. In these representatives the induced action on $Y=\{1,\ldots,s\}\times C$ and its identification with $\Gamma\times_{\Gamma_C}C$ are
\[
 \gamma(j,y)=(k,(t_k^{-1}\gamma t_j)y)\quad
 (\gamma t_j\Gamma_C=t_k\Gamma_C),\qquad
 (j,y)\longmapsto[t_j,y]\in\Gamma\times_{\Gamma_C} C.
\]
Here $\gamma\in\Gamma$, $1\le j\le s$, $y\in C$, and $[th,y]=[t,hy]$ for $t\in\Gamma$, $h\in\Gamma_C$.
\begin{theorem}\label{main}
There exist a connected double quotient $C=\Delta\backslash L/H$ as above and a subgroup $\Gamma_C<\Gamma$ of index $s<\infty$ acting affinely on $C$. On the induced space $Y$ put $\nu=s^{-1}\sum_{j=1}^s\delta_j\otimes m_C$. There are invariant open sets $Y^\circ$ in the smooth locus of $Y$ and $O\subset M$, and an analytic diffeomorphism $i:Y^\circ\to O$, such that
\[
 \nu(Y^\circ)=m(O)=1,\qquad \overline O=M,\qquad
 i_*(\nu|_{Y^\circ})=m,\qquad i(\gamma y)=\gamma i(y)\quad(\gamma\in\Gamma,\ y\in Y^\circ).
\]
\end{theorem}

A longer version, \cite{classification}, classifies all invariant Borel probability measures for real-analytic actions of the same higher-rank lattices on compact manifolds, including singular and nonergodic probabilities. Each ergodic invariant probability is the image of Haar probability under an equivariant analytic immersion from an invariant open subset of full Haar measure in a homogeneous model. The general invariant probabilities are obtained by integrating these measures in analytic families. When volume is not ergodic, Theorem~1.1 of~\cite{classification} gives invariant open sets covering full volume on each of which the action is, up to a finite quotient, the product of a homogeneous action and the identity.

\paragraph{Introduction.}
The examples of Katok and Lewis~\cite{KL}, obtained by modifying algebraic actions near exceptional sets, explain the role of the open set in Theorem~\ref{main}. Benveniste and Fisher~\cite{BenvenisteFisher} showed that such modifications can destroy every global invariant rigid geometric structure.

We began developing the current approach in 2014--2015. From the outset we considered both invariant probabilities for lattice actions and probabilities invariant under parabolic subgroups in the suspension.

This paper is a companion to~\cite{germ}. Its germ superrigidity theorem gives convergent coordinate germs depending measurably on their centers. They express the lattice action through an algebraic representation and a compact cocycle in its centralizer. In Part~\ref{germs} we incorporate the compact cocycle by passing to an extension that records a choice of frame at each center.

For stationary measures, Brown, Eskin, Filip and Rodriguez Hertz~\cite{BEFRH} obtain additional invariance of conditional measures from quantitative nonintegrability. Here we compare the coordinate germs on the suspension using the invariant volume. Intersections of stable manifolds of Cartan elements give the coarse Lyapunov manifolds. By the entropy formulas of Brown, Rodriguez Hertz and Wang~\cite{Entropy}, the conditionals on the vertical coarse Lyapunov plaques are equivalent to plaque volume. For displacements in one coarse Lyapunov space, a single ergodic Cartan element fixes every displacement. The law of each displaced center is invariant under this element. Averaging over the displacements gives a measure absolutely continuous with respect to the suspended volume on the chosen component. Ergodicity and continuity in the displacement then identify each law with that volume. We can therefore pair a chart whose displaced center is $y$ with a chart centered at $y$. On an ergodic component of this joining, every Taylor coefficient of the change of coordinates is constant. Analyticity gives one change of coordinates on that component. Part~\ref{agreement} composes these changes to construct analytic sections of the frame bundle.

In Part~\ref{assembly} the first two derivatives of the charts give a coframe. Its dual fields have constant bracket coefficients and integrate to a Lie group preserving the frame probability. The frame probability is carried by finitely many orbits, each a quotient by a lattice. Dividing by the compact frame group gives the double quotient in Theorem~\ref{main}.

The measure-theoretic and geometric precursors come from nonuniform measure rigidity. Kalinin, Katok and Rodriguez Hertz~\cite{KKRH} obtained absolute continuity from entropy and conditional measures along Lyapunov foliations. Katok and Rodriguez Hertz~\cite{KRHArithmetic} assembled the resulting affine leaf geometry into an arithmetic model. For real-analytic actions of finite-index subgroups of $\mathrm{SL}(n,\mathbb Z)$ on $\mathbb T^n$, $n\ge3$, with the standard action on $H_1(\mathbb T^n;\mathbb Z)$, their earlier theorem~\cite{KRHAnalytic} gives an analytic conjugacy on an invariant open set to the linear action outside finitely many periodic orbits, after passage to finite index, provided there is an ergodic invariant probability whose support is contained in no ball.

In his 1986 ICM address, Zimmer asked whether homogeneous actions, compact isometric actions and nilmanifold automorphisms, together with their modifications, account for ergodic volume-preserving actions in higher rank~\cite[p.~1248]{ZimmerICM}. His report surveyed results under additional hypotheses such as dimension bounds and preservation of a geometric structure. Theorem~\ref{main} gives an open homogeneous classification for real-analytic lattice actions on a set of full volume.

For smooth torus actions of finite-index subgroups of $\mathrm{SL}(n,\mathbb Z)$, Katok, Lewis and Zimmer~\cite{KatokLewisZimmer} obtained affine conjugacy from an Anosov element and an absolutely continuous invariant probability.
In proving local rigidity for algebraic Anosov actions of higher-rank abelian groups, Katok and Spatzier used recurrence along a Lyapunov hyperplane to construct a group of isometries acting transitively on a coarse Lyapunov leaf~\cite[\S2.2.3, Step~4]{KatokSpatzierDifferential}.
Margulis and Qian~\cite{MargulisQian} proved local rigidity for weakly hyperbolic homogeneous and nilmanifold actions. On nilmanifolds, Brown, Rodriguez Hertz and Wang~\cite{BRHW} obtained semiconjugacy to the affine model from a lift and hyperbolic linear data, and smooth conjugacy for Anosov actions. Goetze and Spatzier~\cite{GoetzeSpatzier} classified volume-preserving Cartan actions of higher-rank groups and their irreducible lattices. For semisimple Lie groups with all simple factors of real rank at least two, Damjanovi\'c, Spatzier, Vinhage and Xu~\cite[Theorem~2.5]{DSVX} obtained a model on a double quotient, up to finite cover, under total partial hyperbolicity, accessibility and preservation of volume by a split Cartan subgroup.

Feres~\cite{FeresAffine} identified toral and isometric factors for connection-preserving actions under dimension bounds. For nonisometric actions of lattices in $\mathrm{SL}(n,\mathbb R)$, $n\ge3$, preserving a connection and volume in dimension $\le n+1$, Zeghib~\cite{ZeghibAffine} obtained the standard torus action or its product with the identity on a circle, up to finite cover and finite index. In~\cite{ZeghibAnalytic} he extended a fixed-point linearization to an invariant open set. Feres and Labourie~\cite{FeresLabourie} proved topological cocycle superrigidity for Anosov actions.

\paragraph{Outline of the proof.}
\begin{enumerate}[leftmargin=*,label=\arabic*.]
\item \textbf{Construction and convergence of the coordinate germs.}
Proposition~\ref{prop:group-reduction} passes to a normal finite-index subgroup $\Lambda\lhd\Gamma$, identifies it with a lattice in a semisimple group $G_0$ without compact factors, and fixes one $\Lambda$-ergodic component $m_1$ of the original volume $m$. Put $d=\dim M$. Choose finitely many coordinate charts $\varphi_j:U'_j\to\mathbb R^d$ and open sets $U_j\Subset U'_j$ covering $M$. For $x\in M$ take the least $j$ with $x\in U_j$ and set $H_x(w)=\varphi_j^{-1}(\varphi_j(x)+w)$ near $w=0$. In these coordinates the cocycle is
\[
 T_{\gamma,x}=H_{\gamma x}^{-1}\circ\gamma\circ H_x,\qquad
 T_{\gamma\delta,x}=T_{\gamma,\delta x}\circ T_{\delta,x}
 \quad(\gamma,\delta\in\Lambda).
\]
The companion paper applies Zimmer's cocycle superrigidity~\cite{ZimmerBook} in the form of Fisher--Margulis~\cite[Theorem~1.5]{FisherMargulis} to the finite jets of $T_{\gamma,x}$. It makes these reductions compatible and uses spherical averaging to obtain convergent coordinate germs reducing the cocycle to an algebraic representation times a compact cocycle in its centralizer~\cite[Part~I, Theorem~1.1 and Proposition~2.2]{germ}.

Proposition~\proofplace{action-coordinates} realizes the compact cocycle on an ergodic principal compact extension $c:(P,\eta)\to(M,m_1)$. Its coordinate space is a real vector space $V$ of dimension $d$, the representation is $\rho:G_0\to\GL(V)$ and the compact frame group satisfies $K\le Z_{\GL(V)}(\rho(G_0))$. The resulting germs satisfy
\[
 \psi_p:(V,0)\to(M,c(p)),\qquad
 \gamma\psi_p=\psi_{\gamma p}\rho(\gamma),\qquad
 \psi_{ph}=\psi_p h,\qquad c_*\eta=m_1
 \quad(\gamma\in\Lambda,\ h\in K).
\]
For $\eta$-almost every $p\in P$, the chart $\psi_p$ and its inverse are analytic near $0\in V$ and $c(p)\in M$, respectively. Part~\ref{agreement} determines how these charts vary as the center moves.

\item \textbf{Agreement as the center moves.}
Start with one chart and move its center. Write $\widehat X=(G_0\times P)/\Lambda$ for the suspension of the frame bundle and $X=(G_0\times M)/\Lambda$ for the suspension of $M$. Thus $z=[g,p]=[g\gamma^{-1},\gamma p]$ for $\gamma\in\Lambda$. Its chart is $H_z(w)=[g,\psi_p(\rho(g)^{-1}w)]\in X$ for $w$ near $0$ in $V$. Let $\widehat m$ be the suspension of $\eta$ and $\bar m_1$ that of $m_1$ as in \eqref{susp}. Write $V_\lambda$ for the weight spaces of a maximal split torus. For a nonzero weight $\chi$, fix a small displacement $s\in S_\chi=\bigoplus_{t>0}V_{t\chi}$. For $z$ sampled from $\widehat m$, both $H_z(s)$ and $H_z(0)$ have law $\bar m_1$. This uses equivalence of the conditionals to volume on the coarse Lyapunov plaques and an ergodic Cartan element $b$ fixing $S_\chi$ pointwise (Lemma~\ref{coarse} and Proposition~\ref{eval}).

With $s$ fixed, sample $y\in X$ from $\bar m_1$. Given $y$, independently sample $z\in\widehat X$ conditional on $H_z(s)=y$ and $z'\in\widehat X$ conditional on $H_{z'}(0)=y$ (formula~\eqref{join}, p.~\pageref{join}). For almost every component in the ergodic decomposition of the joining under $(z,z')\mapsto(bz,bz')$, invariance and ergodicity make every Taylor coefficient of $H_z^{-1}\circ H_{z'}$ constant. The same analytic change of coordinates $u$ thus gives $H_{z'}=H_z\circ u$, $u(0)=s$, for almost every pair on that component. Both marginals of each such component are $\widehat m$, and the changes from any two such components differ by right composition with one $h\in K$ (Lemma~\ref{joining}). The group law first appears when we move the center twice. If the second change is $v$, the final chart is $(H_z\circ u)\circ v=H_z\circ(u\circ v)$. Reversing the first move gives $u^{-1}$.

For this Cartan element $b$, the changes commuting with $\rho(b)$ form a local Lie group whose orbit of $0$ contains all small displacements in $S_\chi$ (Lemma~\ref{lie}). Conjugating by $G_0$ supplies motions spanning the nontrivial summands of $V$. In the space fixed by $G_0$, averaging first puts the evaluation law in the finite set of suspended component measures for Lebesgue-almost every displacement. Continuity and connectedness give the same law as at zero throughout a small ball. The calculation in Section~\ref{sec:analytic-sections} gives this finite choice and supplies the remaining motions. Choose one-parameter families of changes whose velocities at $0$ form a basis of $V$. Starting from one frame, their successive times give analytic coordinates for the endpoint by the inverse function theorem. The chart obtained there is therefore analytic in its center and agrees with a frame of Part~\ref{germs} for $m$-almost every center in its open domain (Claim~\ref{cl:section-word}). Two such local choices differ by an element of $K$ depending analytically on the center (Lemma~\ref{sections}).

\item \textbf{Assembly into the homogeneous model.}
Lemma~\ref{finite} first compares these frames with their translates under the stabilizer of $m_1$ in $\Gamma$. After one change of coordinates in $V$, each lift acts on $V$ by a linear map. Two successive lifts and the lift of the product differ by right translation by an element of $K$. Lemma~\ref{coframe} constructs the connection and the coframe $\vartheta$ from the first and second derivatives of $\psi_p$. On the analytic frame bundle $P_1$, its value is an isomorphism
\[
 \vartheta_p:T_pP_1\to V\oplus\mathfrak k,\qquad\mathfrak k=\Lie K.
\]
For $a\in V\oplus\mathfrak k$, let $X_a(p)=\vartheta_p^{-1}a$. Covariance under $\rho(\gamma)\oplus I$ makes the distribution of $(a,b)\mapsto\vartheta_p([X_a,X_b]_p)$ for $p$ sampled from $\eta$ invariant under the induced $\Lambda$-action. Lemma~\ref{H-lem:fixed-law} and ergodicity give a single map for $\eta$-almost every $p\in P_1$. Continuity and full support give constancy on $P_1$. The Jacobi identity for the fields makes this map a Lie bracket.

Part~\ref{assembly} puts the motions from the different directions into one Lie group. A one-parameter change of chart from Part~\ref{agreement} moves each frame with the same coframe velocity (Section~\ref{sec:coframe-from-jets}). In the suspension, higher rank and recurrence extend the local flow relations to a set of full suspended frame measure for each fixed pair of times. Theorem~\proofplace{H-thm:coframe-integration} integrates these relations, and Lemma~\ref{adjoin} includes the compact frame group $K$. On a frame orbit $\Delta\backslash L$, moving the frame is now right multiplication, $\Delta g\mapsto\Delta ga$ for $a\in L$. The stabilizer $\Delta<L$ is discrete, and the finite invariant measure on the orbit makes it a lattice (Lemma~\ref{H-lem:coframe-measured-realization}).

The frame embedding has open image of full Haar measure in the union of these orbits (Lemma~\ref{H-lem:coframe-measured-realization}). Dividing by the free compact frame action gives an analytic diffeomorphism on the corresponding open subset of $\Delta\backslash L/K$. The lattice normalizes the measured group action and hence acts affinely on the model. Finally, the finitely many translates associated to the components $m_1,\ldots,m_k$ recover $\Gamma\curvearrowright(M,m)$. Ergodicity makes the connected model components a finite transitive family. The \hyperlink{proof-main-start}{proof of Theorem~\ref*{main}} identifies their induced action and their equal Haar masses.
\end{enumerate}

\manuscriptpart{The lattice and its coordinate germs}\label{germs}
\noindent\textbf{Outline of Part I.}
Theorem~1.1 in Part~I of~\cite{germ} supplies the analytic coordinates for one ergodic component of a lattice without compact factors. Proposition~\ref{prop:group-reduction} gives that component and records how $\Gamma$ permutes its translates. In Proposition~\ref{action-coordinates} a frame over $x$ specifies which chart at $x$ to use. Changing the frame composes the chart on the right by an element of one compact group. Lemmas~\ref{H-lem:fixed-law} and~\ref{H-lem:hyperplane} give the two ergodicity arguments used to compare these choices.

\section{Reduction to one lattice component}\label{group-reduction}
If $\dim M=0$, then $M$ is finite and ergodicity with full support gives one $\Gamma$-orbit. The theorem follows with $C$ a point. If $\Gamma$ is finite, its invariant Borel sets are the inverse images of the Borel subsets of $M/\Gamma$. The pushforward of $m$ to this compact metrizable quotient takes only the values $0,1$, so it is a point~mass. Its preimage is a closed finite orbit of full measure, hence equals $M$. Theorem~\ref{main} follows with $C$ a point. Assume henceforth that $\dim M>0$ and $\Gamma$ is infinite, and retain the group $G$, lattice $\Gamma$, manifold $M$ and ergodic invariant probability $m$ of Theorem~\ref{main}.

\begin{proposition}[Reduction to one lattice component]\label{prop:group-reduction}
There exist a normal subgroup $\Lambda\lhd\Gamma$ of finite index and a connected semisimple real linear group $G_0$ without compact factors such that $\Lambda$ is isomorphic to a lattice in $G_0$. The group $G_0$ is the connected real group of a simply connected semisimple real algebraic group, and every simple factor has real rank at least two. Under the lattice identification, for each $t\in\Gamma$ there exists $g_t\in G_0$ such that
\[
 t\gamma t^{-1}=g_t\gamma g_t^{-1}\qquad(\gamma\in\Lambda).
\]
There are an integer $1\le k\le[\Gamma:\Lambda]$, distinct $\Lambda$-invariant ergodic probabilities $m_1,\ldots,m_k$, and pairwise disjoint $\Lambda$-invariant Borel sets $X_1,\ldots,X_k\subset M$ satisfying
\begin{equation}\label{components}
 m=\frac1k\sum_{j=1}^k m_j,\qquad
 m_j=k\,m|_{X_j},\qquad m_j(X_j)=1\quad(1\le j\le k).
\end{equation}
The original lattice acts transitively on the component measures,
\[
 \{t_*m_1:t\in\Gamma\}=\{m_1,\ldots,m_k\}.
\]
\end{proposition}
\begin{proof}

The simply connected covering of the semisimple algebraic group with the given Lie algebra induces a finite central covering $p:\widehat G=G_0\times C_0\to G$ on the connected real groups, since $G$ is linear~\cite{BorelAlgebraic}. Here $C_0$ consists of the compact factors, and $G_0$ is the connected real group from the noncompact factors.
Let $\pi:\widehat G\to G_0$ be the projection, $\widehat\Gamma=p^{-1}(\Gamma)$ and $\Delta=\pi(\widehat\Gamma)$. If distinct $g_j\in\Delta$ converge to $1$, choose $(g_j,c_j)\in\widehat\Gamma$ and a subsequence $c_{j_m}\to c$ in $C_0$. Then
\[
 1\ne (g_{j_m},c_{j_m})^{-1}(g_{j_{m+1}},c_{j_{m+1}})\to1,
\]
contradicting discreteness. Thus $\Delta$ is discrete and $\ker(\pi|_{\widehat\Gamma})=\widehat\Gamma\cap(\{1\}\times C_0)$ is finite. The quotient map $\widehat\Gamma\backslash\widehat G\to\Delta\backslash G_0$ sends invariant probability to invariant probability, which gives finite~covolume for $\Delta$.

The lattice $\widehat\Gamma$ is finitely generated~\cite[Remark~6.18]{Raghunathan}. By Selberg's lemma~\cite{Selberg} choose $\Lambda^{\prime}<\widehat\Gamma$ torsion-free of finite index. The subgroup $\widehat\Lambda=\bigcap_{h\in\widehat\Gamma}h\Lambda^{\prime}h^{-1}$ is normal and has finite index since $h\Lambda^{\prime}h^{-1}$ depends only on $h\Lambda^{\prime}$. For both the covering map $p:\widehat G\to G$ and $\pi$ we have
\[
 \ker(p|_{\widehat\Lambda})=\ker(\pi|_{\widehat\Lambda})=\{1\},\qquad
 \Lambda=p(\widehat\Lambda)\triangleleft\Gamma,\qquad
 [\Gamma:\Lambda]\le[\widehat\Gamma:\widehat\Lambda]<\infty .
\]
Since $[\Delta:\pi(\widehat\Lambda)]\le[\widehat\Gamma:\widehat\Lambda]<\infty$, $\pi(\widehat\Lambda)$ is a lattice in $G_0$. The isomorphism $\iota:\Lambda\to\pi(\widehat\Lambda)$ is given by $\iota(p(\hat\gamma))=\pi(\hat\gamma)$ for $\hat\gamma\in\widehat\Lambda$. If $p(g_t,c_t)=t$, then for $\gamma=p(\hat\gamma)\in\Lambda$,
\[
 \iota(t\gamma t^{-1})
 =\pi\bigl((g_t,c_t)\hat\gamma(g_t,c_t)^{-1}\bigr)
 =g_t\iota(\gamma)g_t^{-1}.
\]
We henceforth use this identification.
Write $m=\int\xi\,d\tau(\xi)$ for the ergodic decomposition under $\Lambda$ \cite[Theorems~4.2 and~4.4]{Varadarajan}. Since $\Lambda\lhd\Gamma$, the finite group $\Gamma/\Lambda$ acts on its component measures by $t\cdot\xi=t_*\xi$. The component map is equivariant modulo $m$, by uniqueness of disintegration. Thus $\tau$ is ergodic. Its image on the standard Borel space of finite orbits is a point~mass, since every set in that quotient is invariant. Consequently
\[
 \tau=\frac1k\sum_{j=1}^k\delta_{m_j},\qquad
 \{m_1,\ldots,m_k\}=\{t_*m_1:t\in\Gamma\}.
\] Choose disjoint $\Lambda$-invariant Borel sets $X_j$ of full $m_j$-measure. Invariance under $\Gamma$ gives
\[
 m=\frac1k\sum_{j=1}^k m_j,\qquad m_j=k\,m|_{X_j}.
\]
The orbit has at most $[\Gamma:\Lambda]$ members, since $\Lambda$ fixes each component. This proves all the assertions.
\end{proof}
For Parts~\ref{germs} and~\ref{agreement}, fix the component $m_1$ and the identification $\Lambda<G_0$ of Proposition~\ref{prop:group-reduction}. Lemma~\ref{sections} will realize $m_1$ as normalized volume on an invariant open subset of $M$.

\section{Analytic germs and compact frames}

\begin{definition}[Coordinate germ]\label{def:coordinate-germ}
A coordinate germ $\psi:(V,0)\to(M,x)$ is the germ at $0$ of an analytic local diffeomorphism from a neighborhood of $0$ in a real vector space $V$ onto a neighborhood of $x$, with $\psi(0)=x$. Two representatives give the same germ if they agree on an open neighborhood of $0$. A family $p\mapsto\psi_p$ is measurable if its centers $p\mapsto\psi_p(0)$ are measurable and, in each target chart, $p\mapsto D_0^j\psi_p$ is measurable for every $j\ge1$ on the set of centers in that chart.
\end{definition}
\begin{definition}[Principal compact extension]\label{def:compact-extension}
A principal $K$-extension of $\Lambda\curvearrowright(M,m_1)$ consists of a probability-preserving $\Lambda$-action on $(P,\eta)$, a measurable map $c:P\to M$ with $c_*\eta=m_1$, and a free right action of a compact group $K$ preserving $\eta$. Over a $\Lambda$-invariant set of full $m_1$-measure,
\[
 c^{-1}(c(p))=pK,\qquad c(\gamma p)=\gamma c(p),\qquad
 \gamma(ph)=(\gamma p)h \quad(\gamma\in\Lambda,\ h\in K).
\]
We write $m_K$ for Haar probability on $K$.
\end{definition}

For a subspace $S$ of a finite-dimensional inner-product space put $B_r(S)=\{v\in S:\|v\|<r\}$ for $r>0$. In the next proposition the norm on $V$ is chosen invariant under the compact frame group.
\begin{proposition}[The ergodic compact frame extension]\label{action-coordinates}
There are a real vector space $V$ of dimension $\dim M$, an algebraic representation
$\rho:G_0\to\mathrm{GL}(V)$, a compact subgroup
$K<\mathrm{GL}(V)$ centralizing $\rho(G_0)$, and an ergodic
probability-preserving principal $K$-extension
$c:(P,\eta)\to(M,{m_1})$ with measurable analytic diffeomorphism germs
$\psi_p:(V,0)\to(M,c(p))$. On a $\Lambda$-invariant set of full $m_1$-measure, for every frame $p$ over that set, $\gamma\in\Lambda$, and $h\in K$,
\begin{equation}\label{frames}
 \begin{gathered}
 \gamma\psi_p=\psi_{\gamma p}\rho(\gamma),\qquad
 \psi_{ph}=\psi_p\circ h,\qquad
 \gamma(ph)=(\gamma p)h,\\
 c_*\eta=m_1,\qquad \gamma_*\eta=\eta,\qquad
 (R_h)_*\eta=\eta,\qquad R_h(p)=ph .
 \end{gathered}
\end{equation}
For $m_1$-almost every $x\in M$ and $p,p'\in c^{-1}(x)$, the conditional probability $\eta_x$ in $\eta=\int_M\eta_x\,dm_1(x)$ and the derivatives $D_0\psi_p:V\to T_xM$ satisfy
\[
  \eta_x=(h\mapsto ph)_*m_K,\qquad D_0\psi_p=D_0\psi_{p'}\ \Longrightarrow\ p=p'.
\]
For $m_1$-almost every $x\in M$ there is a measurable $r_x>0$ such that all $\psi_p$, $c(p)=x$, embed $B_{r_x}(V)$ onto the same neighborhood of $x$.
\end{proposition}
\begin{proof}
Apply~\cite[Part~I, Theorem~1.1]{germ} to the compact analytic $M$, the reduced lattice and the ergodic invariant probability ${m_1}$ of Proposition~\ref{prop:group-reduction}. It gives $V=\mathbb R^{\dim M}$, a real algebraic $\rho:G_0\to\GL(V)$, a compact $K_0\le Z_{\GL(V)}(\rho(G_0))$, a measurable cocycle $\kappa:\Lambda\times M\to K_0$ with $\kappa(\gamma\delta,x)=\kappa(\gamma,\delta x)\kappa(\delta,x)$ for $\gamma,\delta\in\Lambda$, and analytic coordinate germs $\Phi_x:(V,0)\to(M,x)$ with
\begin{equation}\label{frame-first-order}
 \gamma\Phi_x=\Phi_{\gamma x}\rho(\gamma)\kappa(\gamma,x),\qquad
 D_x\gamma\,D_0\Phi_x=D_0\Phi_{\gamma x}\rho(\gamma)\kappa(\gamma,x).
\end{equation}
The identities hold for all $\gamma\in\Lambda$ on one invariant set of full ${m_1}$-measure. Average an inner product on $V$ over $K_0$ and use its norm in $B_r(S)$.
\par\Needspace{3\baselineskip}\smallskip\noindent\textbf{Reduction of the compact group.} Use the action and probability on $M\times K_0$ given by
\[
 \gamma(x,h)=(\gamma x,\kappa(\gamma,x)h),\qquad m_1\otimes m_{K_0},
\]
 where $m_{K_0}$ is normalized Haar measure.
Let $\chi$ be the ergodic decomposition law on the space $\mathscr E_\Lambda(M\times K_0)$ of ergodic invariant probabilities~\cite[Theorems~4.2 and~4.4]{Varadarajan}, so ${m_1}\otimes m_{K_0}=\int\lambda\,d\chi(\lambda)$. Since ${m_1}=\int(\operatorname{pr}_M)_*\lambda\,d\chi(\lambda)$ is extreme among invariant probabilities, $(\operatorname{pr}_M)_*\lambda={m_1}$ for $\chi$-almost every $\lambda\in\mathscr E_\Lambda(M\times K_0)$. Fix one such probability and denote it by $\eta_0$. For every $\gamma\in\Lambda$, uniqueness of disintegration gives, for $m_1$-almost every $x\in M$,
\[
 \eta_0=\int\eta_{0,x}\,d{m_1}(x),\qquad
 \eta_{0,\gamma x}=\kappa(\gamma,x)_*\eta_{0,x}.
\]
Here $h_*\lambda$ denotes the image of $\lambda$ under $k\mapsto hk$. The orbit space of this $K_0$-action on $\operatorname{Prob}(K_0)$ is compact and metrizable. Since $K_0\eta_{0,\gamma x}=K_0\kappa(\gamma,x)_*\eta_{0,x}=K_0\eta_{0,x}$, ergodicity gives $\eta_*\in\operatorname{Prob}(K_0)$ with ${m_1}\{x\in M:K_0\eta_{0,x}=K_0\eta_*\}=1$. Choose measurable $b:M\to K_0$ with $\eta_{0,x}=b(x)_*\eta_*$ there.
Indeed the map $h\mapsto h_*\eta_*$ is a continuous surjection from $K_0$ to its orbit. Its graph has compact fibers, so it has a Borel section by~\cite[Theorem~18.18]{Kechris}.
Define
\[
 K=\{k\in K_0:k_*\eta_*=\eta_*\},\qquad
 \kappa_K(\gamma,x)=b(\gamma x)^{-1}\kappa(\gamma,x)b(x)\in K .
\]
On $M\times K_0$ the Borel bijection $B(x,k)=(x,b(x)^{-1}k)$ has inverse $B^{-1}(x,k)=(x,b(x)k)$. Then $B_*\eta_0={m_1}\otimes\eta_*$ and $B\gamma B^{-1}(x,k)=(\gamma x,\kappa_K(\gamma,x)k)$. Since $\kappa_K(\gamma,x)\in K$, the coset satisfies $K\kappa_K(\gamma,x)k=Kk$.
Ergodicity gives $k_0\in K_0$ with $\eta_*(Kk_0)=1$. The defining identity $k_*\eta_*=\eta_*$ for all $k\in K$ gives $\eta_*=(k\mapsto kk_0)_*m_K$.
Right translation by $k_0^{-1}$ identifies the conditional measure on $Kk_0$ with $m_K$. The resulting ergodic action is
\[
 P=M_*\times K,\quad \eta={m_1}\times m_K,\quad c(x,k)=x,
 \qquad \gamma(x,k)=(\gamma x,\kappa_K(\gamma,x)k),
\]
where $M_*=\bigcap_{\gamma\in\Lambda}\gamma M_1$. Here $M_1$ is a Borel set of ${m_1}$-measure one on which the preceding identities hold for every $\gamma\in\Lambda$. Countability gives ${m_1}(M_*)=1$ and $\gamma M_*=M_*$.

Set
\begin{equation}\label{actual-frame-charts}
 \psi_{(x,k)}=\Phi_x\circ b(x)k .
\end{equation}
Substitute $\kappa(\gamma,x)b(x)=b(\gamma x)\kappa_K(\gamma,x)$ into \eqref{frame-first-order}. Since $\rho(G_0)$ commutes with $K_0$,
\[
 \gamma\psi_{(x,k)}
 =\Phi_{\gamma x}\rho(\gamma)\kappa(\gamma,x)b(x)k
 =\Phi_{\gamma x} b(\gamma x)\kappa_K(\gamma,x)k\,\rho(\gamma)
 =\psi_{\gamma(x,k)}\rho(\gamma).
\]
Right translation preserves $m_K$ and $\psi_{(x,kh)}=\psi_{(x,k)}h$ for $h\in K$. The fiber coordinate is determined by the derivative through
\[
 k=b(x)^{-1}(D_0\Phi_x)^{-1}D_0\psi_{(x,k)} .
\]
By~\cite[Part~I, Theorem~1.1]{germ}, choose a measurable $r_x>0$ on which $\Phi_x$ is an analytic embedding. Indeed the theorem gives measurable radii for the germ and its inverse, and one can shrink the source radius until its image is in the inverse domain. Orthogonality gives $b(x)kB_{r_x}(V)=B_{r_x}(V)$ for every $k\in K$, so $\psi_{(x,k)}(B_{r_x}(V))=\Phi_x(B_{r_x}(V))$.
\end{proof}

\section{Invariant probabilities in a representation}

\begin{lemma}[Invariant probabilities in a representation]\label{H-lem:fixed-law}
Let $G$ be connected, semisimple and real linear with finite center and
no compact factors, $\Gamma<G$ a lattice, and $W$ a finite-dimensional
algebraic $G$-module. Fix $0\le j\le\dim W$. For $Y=W$ or $\Gr_j(W)$ put
$Y^G=\{y\in Y:gy=y\ (g\in G)\}$. For a Borel probability $\lambda$ on $Y$, $\gamma_*\lambda=\lambda$ for all $\gamma\in\Gamma$ implies $\lambda(Y^G)=1$. Let $\Gamma$ preserve an ergodic probability $\mu$ on $X$. If $s:X\to Y$ is measurable and $s(\gamma x)=\gamma s(x)$ for each $\gamma\in\Gamma$ and $\mu$-almost every $x\in X$, then $s_*\mu=\delta_y$ for some $y\in Y^G$.
\end{lemma}
\begin{proof}
Let $m_{G/\Gamma}$ be the $G$-invariant probability on $G/\Gamma$. Since $\gamma_*\lambda=\lambda$, the probability
\[
 \bar\lambda=\int_{G/\Gamma}g_*\lambda\,dm_{G/\Gamma}(g\Gamma)
\]
is well defined and $G$-invariant. Embed $W$ in $\mathbb P(W\oplus\mathbb R)$ by $w\mapsto[w:1]$, with $G$ trivial on the last coordinate. For the Grassmannian use the Pl\"ucker embedding $S\mapsto[\bigwedge^jS]$ (for $j=0$ the assertion is immediate). Both maps are equivariant and injective. The restricted-root unipotent groups generate $G$. Finitely many one-parameter groups suffice, by taking bases in the root spaces. Let $\exp N$ be the time-one map of one such subgroup on the projective representation, and choose a nonzero vector $w$ in its representation space ($W\oplus\mathbb R$ or $\bigwedge^j W$). With $a=\max\{i\in\mathbb Z_{\ge0}:N^iw\ne0\}$,
\[
 [\exp(nN)w]=\left[\sum_{i=0}^a\frac{n^i}{i!}N^i w\right]\longrightarrow[N^a w].
\]
At a recurrent line this limit is $[w]$. If $a>0$, then $w$ is an eigenvector of the nilpotent operator $N^a$ with nonzero eigenvalue. Thus $a=0$ and $Nw=0$. Poincar\'e recurrence for $\bar\lambda$ and intersection over the chosen generators give
\[
 1=\bar\lambda(Y^G)=\int_{G/\Gamma}\lambda(g^{-1}Y^G)\,dm_{G/\Gamma}(g\Gamma)=\lambda(Y^G).
\]
For $\lambda=s_*\mu$, equivariance gives $s(\gamma x)=s(x)$ for $\mu$-almost every $x\in X$ and each $\gamma$. Ergodicity gives $s_*\mu=\delta_y$ for some $y\in Y^G$.
\end{proof}

\subsection*{Ergodic elements in a split torus}
Let a connected semisimple group $G$ with finite center and without compact factors act ergodically on a probability space $(X,\mu_X)$. Assume $g_*\mu_X=\mu_X$ for $g\in G$ and strong continuity of the Koopman representation $U_g F=F\circ g^{-1}$ on $L^2(\mu_X)$. Write $\mathfrak a$ for the Lie algebra of a maximal real split torus. Assume $\dim\mathfrak a_i=\operatorname{rank}_{\mathbb R}G_i\ge2$ for each simple factor $G_i$, where $\mathfrak a_i$ is the projection of $\mathfrak a$ on its Lie algebra.

\begin{lemma}[Ergodicity on a weight hyperplane]\label{H-lem:hyperplane}
If $H\in\mathfrak a$ projects nontrivially to every simple factor, the $\exp H$-invariant functions in $L^2(X,\mu_X)$ are the constants, $\ker(U_{\exp H}-I)=\mathbb C\mathbf1$. For each $\chi\in\mathfrak a^*$ there exists such an $H\in\ker\chi$.
\end{lemma}
\begin{proof}
Fix $H\in\mathfrak a=\bigoplus_j\mathfrak a_j$ with nonzero projection to each simple factor. For $a=\exp H$, $Y\in\operatorname{Lie}(G)$ and $[H,Y]=\lambda Y$ with $\lambda\in\mathbb R$, we have
\[
 a^{-n}\exp(tY)a^n=\exp(e^{-n\lambda}tY),
 \qquad t\in\mathbb R,\quad n\in\mathbb Z.
\]
For $\lambda\ne0$ and $v\in L^2(X,\mu_X)$ with $U_av=v$, unitarity gives $\|U_{\exp(tY)}v-v\|_2=\|U_{\exp(e^{-n\lambda}tY)}v-v\|_2\to0$ when $n\lambda\to+\infty$, by strong continuity. In each simple factor, the
subalgebra generated by the nonzero $H$-weight spaces is a
nonzero ideal because the zero weight space normalizes every nonzero weight space and all nonzero weight spaces belong to the subalgebra. It is therefore the whole factor. Thus $U_gv=v$ for every $g\in G$, so $v\in\mathbb C\mathbf1$ by ergodicity. Finally, given $\chi\in\mathfrak a^*$, each $\ker\operatorname{pr}_j$ has codimension $\dim\mathfrak a_j\ge2$ in $\mathfrak a$, whereas $\ker\chi$ has codimension at most one. The finitely many proper subspaces $\ker\chi\cap\ker\operatorname{pr}_j$ cannot cover $\ker\chi$. Any $H$ outside their union in $\ker\chi$ has the required projections.
\end{proof}

\manuscriptpart{Compatibility of the coordinate germs by two-point joinings}\label{agreement}
\section{Analytic sections of the compact frame extension}
For the reduced lattice $\Lambda<G_0$ of Proposition~\ref{prop:group-reduction}, the ergodic decomposition \eqref{components} of the original volume $m$ is
\[
 m=\frac1k\sum_{j=1}^k m_j,\qquad m_j=k\,m|_{X_j}\quad(1\le j\le k),
\]
where the $X_j\subset M$ are pairwise disjoint $\Lambda$-invariant Borel sets. Fix $m_1$. Proposition~\ref{action-coordinates} gives an ergodic principal $K$-extension $c:(P,\eta)\to(M,m_1)$ and a representation $\rho:G_0\to\mathrm{GL}(V)$, with $\dim V=\dim M$ and $K$ compact in the centralizer of $\rho(G_0)$. The analytic coordinate germs of \eqref{frames} satisfy
\[
 \psi_p:(V,0)\longrightarrow (M,c(p)),\qquad
 \psi_{\gamma p}=\gamma\psi_p\rho(\gamma)^{-1},\qquad \psi_{ph}=\psi_ph
 \quad(\gamma\in\Lambda,\ h\in K).
\]
The center $c(p)$ and the derivative $D_0\psi_p$ determine the frame $p$. For $m_1$-almost every $x\in M$ and any $p\in c^{-1}(x)$, the conditional probability of $\eta$ on $c^{-1}(x)$ is $(h\mapsto ph)_*m_K$, where $m_K$ is Haar probability on $K$.

\begin{lemma}[Analytic sections of the frame bundle]\label{sections}
There is a $\Lambda$-invariant open set $O_1\subset M$ with $m_1(O_1)=1$ and $m_1|_{O_1}=k\,m|_{O_1}$. Each point of $O_1$ has an open neighborhood $U\subset O_1$ and an analytic map $(x,w)\mapsto s_x(w)$ defined near $U\times\{0\}$ such that $s_x(0)=x$ and
\[
 \{\,\psi_p:p\in c^{-1}(x)\,\}=\{\,s_xh:h\in K\,\}
 \qquad\text{for }m_1\text{-almost every }x\in U .
\]
The sections define an analytic principal $K$-bundle over $O_1$. Over a set of full $m_1$-measure the identification $p\mapsto\psi_p$ intertwines both the $\Lambda$- and $K$-actions. In the coordinates $(x,h)\mapsto s_x\circ h$ the frame measure is $dm_1(x)\,dm_K(h)$.
\end{lemma}

\noindent\textbf{Proof outline.}
The proof of Lemma~\ref{sections}, the main result of this part, is completed in Section~\proofplace{sec:analytic-sections}. It uses two results proved below. To pair a chart evaluated at a displacement $s$ with a chart centered at the same point, keeping the prescribed frame measure on both coordinates, the two evaluations must have the same law. For $s$ in a coarse Lyapunov space, Lemma~\ref{coarse} identifies the conditional measures on the plaques, and Proposition~\ref{eval} uses them to prove equality of the two laws. With this common law we form the joining~\eqref{join}. On an ergodic component, Lemma~\ref{joining} gives a single change of coordinates between the paired charts. Lemma~\ref{lie} then gives one-parameter families of these changes. In Section~\ref{sec:analytic-sections}, Claim~\ref{cl:fixed-evaluation} treats the directions fixed by the representation. Claim~\ref{cl:section-word} identifies the chart obtained by composing the changes. The center of this chart depends analytically on the parameters, and the center velocities form a basis. Inverting the map from parameters to centers gives $s_x(w)$.

\section{Coarse Lyapunov plaques and evaluation}

Let $X=(G_0\times M)/\Lambda$ and $\widehat X=(G_0\times P)/\Lambda$, where
$(g,p)\sim(g\gamma^{-1},\gamma p)$ for $\gamma\in\Lambda$, and likewise on $G_0\times M$. The base projections are $[g,x]\mapsto g\Lambda$ on $X$ and $[g,p]\mapsto g\Lambda$ on $\widehat X$. For fixed $g\in G_0$ let $\iota_g:M\to X$, $\iota_g(x)=[g,x]$, be the identification of $M$ with the fiber over $g\Lambda$. The suspended probability $\widehat m$ has Haar probability on $G_0/\Lambda$ as marginal and $([g,\cdot])_*\eta$ as its conditional at $g\Lambda$. Define $\bar m_j$ and $\bar m$ in the same way using $m_j$ and $m$. These conditionals are independent of the representative $g$ by $\Lambda$-invariance.

For $z=[g,p]$ define the chart with values in that fiber by
\begin{equation}\label{susp}
 H_z=\iota_g\circ\psi_p\circ\rho(g)^{-1},\qquad
 e_0(z)=H_z(0),\qquad H_{az}=a\circ H_z\circ\rho(a)^{-1}\quad(a\in G_0).
\end{equation}
Indeed $\psi_{\gamma p}=\gamma\psi_p\rho(\gamma)^{-1}$ by \eqref{frames}, so replacing $(g,p)$ by $(g\gamma^{-1},\gamma p)$ gives the same germ. Its center has law $(e_0)_*\widehat m=\bar m_1$ since $c_*\eta=m_1$. A $G_0$-invariant measurable function on $\widehat X$ lifts to a function constant in $g$ for Haar-almost every $g\in G_0$ and $\Lambda$-invariant in $p$ modulo $\eta$. Ergodicity of $\eta$ makes it constant. The same argument gives $G_0$-ergodicity of each $\bar m_j$.

Fix a maximal real split torus with Lie algebra $\mathfrak a$. Write $V=\bigoplus_\lambda V_\lambda$, where $V_\lambda=\{v\in V:\rho(\exp H)v=e^{\lambda(H)}v\ (H\in\mathfrak a)\}$ for each weight $\lambda\in\mathfrak a^*$. For a nonzero weight $\chi$ and $r>0$, put
\[
 S_\chi=\bigoplus_{\lambda=t\chi,\ t>0}V_\lambda\ \subset V,\qquad B_r(S_\chi)=\{s\in S_\chi:\|s\|<r\},
\]
using the norm fixed in Proposition~\ref{action-coordinates}. Choose $H_0\in\ker\chi$ with nonzero projection to the Lie algebra of each simple factor and with $\lambda(H_0)\ne0$ for every weight not proportional to $\chi$. The excluded sets are proper subspaces of $\ker\chi$. Indeed, the kernel of projection onto the split Cartan of any single factor has codimension at least two in $\mathfrak a$. Put $b=\exp H_0$ and $C_b=\Fix(\rho(b))$. In particular $S_\chi\subset C_b$.

\begin{lemma}[Common domain along the fixed space]\label{wall}
Let $b=\exp H_0$ with $H_0$ having nonzero projection to every simple factor, and put $C_b=\Fix(\rho(b))$. Left translation by $b$ is ergodic on $(\widehat X,\widehat m)$ and on each $(X,\bar m_j)$. For some $r>0$ and $\widehat m$-almost every $z\in\widehat X$, $H_z$ has an analytic embedding representative on an open neighborhood in $V$ of $\overline B_{2r}(C_b)$. For fixed $s\in B_{2r}(C_b)$, the evaluation $e_s:\widehat X\to X$, $e_s(z)=H_z(s)$, satisfies $e_s(bz)=b e_s(z)$ for $\widehat m$-almost every $z\in\widehat X$.
\end{lemma}
\begin{proof}
The suspended probabilities are $G_0$-ergodic by the paragraph following \eqref{susp}. Their Koopman operators are strongly continuous. For the measurable frame suspension, this follows first for functions $\varphi(g)f(p)$ in a local trivialization, with $\varphi$ continuous and compactly supported in a local section of $G_0/\Lambda$ and $f$ bounded and measurable, because sufficiently small left translations change only the first factor there. Finite sums of such functions are dense in $L^2$. Lemma~\ref{H-lem:hyperplane} applies to $b$ for each of the stated probabilities.

Choose $r>0$ so that the frames admitting an embedding on $B_{3r}(V)$ have positive measure. By ergodicity, for $\widehat m$-almost every $z\in\widehat X$ there is an $n\ge0$ with $b^n z$ in that set. Then $b^{-n}H_{b^nz}\rho(b)^n$ is an embedding on $\rho(b)^{-n}B_{3r}(V)$, which contains $\overline B_{2r}(C_b)$. Two choices agree as germs at zero by \eqref{susp}, hence throughout the convex intersection of their domains. The first admissible $n\geq0$ is a measurable choice. The same identity with centers $z,bz$ gives the assertion about $e_s$.
\end{proof}

Let $\cV$ be the foliation of $X$ by the fibers of $X\to G_0/\Lambda$. For the chosen nonzero weight $\chi$ we specify the Cartan elements used to define its plaques. Recall that $\chi(H_0)=0$ and $\lambda(H_0)\ne0$ for weights $\lambda$ not proportional to $\chi$. Choose $H_1\in\mathfrak a$ with $\chi(H_1)<0$ and put $a_\pm=\exp(\pm H_0+\epsilon H_1)$. We may take $H_0$ off the restricted-root hyperplanes that do not contain $\ker\chi$, and $H_1$ off the remaining hyperplanes. For small $\epsilon>0$ both $a_\pm$ are regular and have nonzero projection in every simple factor. For weights not proportional to $\chi$ take $|\epsilon\lambda(H_1)|<|\lambda(H_0)|$. Consequently, for every weight $\lambda$ of $V$,
\begin{equation}\label{eq:coarse-sign}
 \lambda(\log a_+)<0\ \text{and}\ \lambda(\log a_-)<0
 \quad\Longleftrightarrow\quad \lambda=t\chi\text{ for some }t>0 .
\end{equation}
Complete $\log a_+,\log a_-$ to a basis of $\mathfrak a$ and let $A_{\mathbb Z}$ be the exponential of its integer span. These vectors are independent because $\chi(H_0)=0$, $H_0\ne0$ and $\chi(H_1)<0$. For $a\in A_{\mathbb Z}$ let $W^s_{\cV,a}(y)$ denote its stable manifold inside the fiber through $y\in X$. For $\bar m_1$-almost every $y\in X$, $W^{[\chi]}_{\cV}(y)$ denotes the local branch through $y$ of
\[
 \bigcap_{\substack{a\in A_{\mathbb Z}\\ \chi(\log a)<0}} W^s_{\cV,a}(y).
\]
The intersection uses a countable group. The coarse Lyapunov manifold and its tangent are given by~\cite[Proposition~4.6, Lemma~4.8 and Definition~4.10]{Entropy} and~\cite[\S3.2]{Boundary}.
For a subordinate measurable partition $\zeta$, the atom $\zeta(y)\subset W^{[\chi]}_{\cV}(y)$ contains a relative neighborhood of $y$ for $\bar m_1$-almost every $y\in X$. Write $(\bar m_1)^\zeta_y$ for its conditional probability on this atom and $\vol_{\zeta(y)}$ for the restriction of plaque volume.

\begin{lemma}[The vertical coarse Lyapunov plaques]\label{coarse}
For each such partition $\zeta$ and each nonzero weight $\chi$,
\[
 \begin{aligned}
 (\bar m_1)^\zeta_y&\sim\vol_{\zeta(y)}&&(\bar m_1\text{-a.e. }y\in X),\\
 \operatorname{germ}_{e_0(z)}W^{[\chi]}_{\cV}(e_0(z))
   &=\operatorname{germ}_{e_0(z)}H_z(S_\chi)&&(\widehat m\text{-a.e. }z\in\widehat X).
 \end{aligned}
\]
\end{lemma}
\begin{proof}
\textbf{Entropy on the plaques.} In the product charts of~\cite[\S2.3]{Boundary}, the fiber plaques of $\cV$ are the slices on which the base coordinate is constant. Their graph over the vertical tangent space is the zero map, so $\cV$ satisfies the tameness condition of~\cite[Def.~4.4]{Entropy}. In the admissible charts of~\cite[Prop.~2.5]{Boundary}, the Cartan derivative on $T\cV$ has weight spaces $D_0H_z V_\lambda$. Indeed, for a regular $a=\exp H$, differentiate \eqref{susp} to get
\[
 D_{e_0(z)}a^n\,D_0H_z v=D_0H_{a^nz}\rho(a)^n v .
\]
Choose $C$ so that
\[
 D=\{z\in\widehat X:\max(\|D_0H_z\|,\|(D_0H_z)^{-1}\|)\le C\},\qquad \widehat m(D)>0 .
\]
The vertical derivative cocycle is over $(X,\bar m_1)$. Its Oseledets regular set pulls back under $e_0$ to a set of full $\widehat m$-measure, since $(e_0)_*\widehat m=\bar m_1$. Choose $a$ with nonzero projection in every simple factor. By Lemma~\ref{H-lem:hyperplane}, for $\widehat m$-almost every $z\in\widehat X$ in this pullback there are returns $a^nz\in D$ as $n\to+\infty$ and as $n\to-\infty$. For $0\ne v\in V_\lambda$ and every return $a^nz$ to this block,
\[
 C^{-1}e^{n\lambda(H)}\|v\|\le\|D_{e_0(z)}a^n\,D_0H_zv\|
          \le C e^{n\lambda(H)}\|v\|.
\]
The limits given by Oseledets' theorem along $n\to+\infty$ and $n\to-\infty$ therefore put $D_0H_zv$ in the $\lambda(H)$-space. A finite set of regular $H$ separating the weights identifies the joint spaces. Proposition~2.5 of~\cite{Boundary} supplies the logarithmic moment and chart estimates also for a nonuniform lattice, since each $\bar m_j$ projects to Haar probability.

Choose $H\in\log A_{\mathbb Z}$ with $\chi(H)>0$ off the root and nonzero weight hyperplanes and with nonzero projection to each simple factor. Such lattice vectors exist since these conditions define a nonempty open cone. Put $a=\exp H$. On every positive coarse class $[\xi]=\{t\xi:t>0\}$ occurring vertically put
\[
 d_\xi=\sum_{\lambda\in[\xi]}\lambda(H)\dim V_\lambda,
 \qquad h_\xi=h_{\bar m_1}(a\mid W^{[\xi]}_{\cV}).
\]
Here $W^{[\xi]}_{\cV}$ is the coarse Lyapunov foliation within the fibers, as in~\cite[\S3.2]{Boundary}. The $A_{\mathbb Z}$-action on $(X,\bar m_1)$ is ergodic because its element $a$ is ergodic by Lemma~\ref{H-lem:hyperplane}. In~\cite[Theorem~13.1]{Entropy} take the foliation $\cV$ by fibers of $X\to G_0/\Lambda$ and the partition with single atom $X$. The intersection for a class $[\xi]$ is $W^{[\xi]}_{\cV}$. If no vertical weight belongs to $[\xi]$, Proposition~4.6 of~\cite{Entropy} gives $T_yW^{[\xi]}_{\cV}(y)=0$ for $\bar m_1$-almost every $y$. The local plaques are points and their entropy contribution is zero. Thus only the vertical classes enter the first sum below. The second equality is the relative Pesin formula.
\[
 \sum_{[\xi]:\,\xi(H)>0}h_\xi
 =h_{\bar m_1}(a\mid\cV)
 =\sum_{\lambda(H)>0}\lambda(H)\dim V_\lambda
 =\sum_{[\xi]:\,\xi(H)>0}d_\xi .
\]
For the relative formula, let $du$ be Haar probability on $G_0/\Lambda$ and take the fundamental-domain section $\sigma:G_0/\Lambda\to G_0$ of~\cite[\S2]{Boundary}. Over $u\mapsto au$ the map of the compact fiber $M$ is the action of $\gamma_a(u)=\sigma(au)^{-1}a\sigma(u)\in\Lambda$. Its sample probability is $m_1\ll\vol_M$ at Haar-almost every $u\in G_0/\Lambda$, where $\vol_M$ denotes Riemannian volume. If $|\cdot|$ is word length in a fixed finite symmetric generating set of $\Lambda$, the chain rule on compact $M$ gives
\[
 \log^+\|\gamma_a(u)^{\pm1}\|_{C^2}\le C(1+|\gamma_a(u)|),
 \qquad \int_{G_0/\Lambda}|\gamma_a(u)|\,du<\infty .
\]
The integral is the return word estimate in~\cite[proof of Prop.~2.5]{Boundary}. The Pesin entropy formula for stationary sequences of maps of a compact manifold with absolutely continuous sample probabilities now gives the vertical sum~\cite{PesinRandom}. Indeed, for a finite partition $\mathcal P$ of $M$ put $R_u^0=\id$ and $R_u^i=\gamma_a(a^{i-1}u)\cdots\gamma_a(u)$ for $i\ge1$. In the trivialization $G_0/\Lambda\times M$, the first $n$ iterates of $\mathcal P$ on the fiber over $u$ give the partition $\bigvee_{i=0}^{n-1}(R_u^i)^{-1}\mathcal P$. Thus the relative entropy is
\[
 h_{\bar m_1}(a\mid\cV)=
 \sup_{\mathcal P}\lim_{n\to\infty}\frac1n
 \int_{G_0/\Lambda} H_{m_1}\!\left(\bigvee_{i=0}^{n-1}(R_u^i)^{-1}\mathcal P\right)\,du .
\]
Here $H_{m_1}(\mathcal Q)=-\sum_{A\in\mathcal Q}m_1(A)\log m_1(A)$, with $0\log0=0$. One may take the supremum over a refining sequence of finite partitions generating the Borel sets of $M$. The integrand depends only on the return sequence, whose sample measure is $m_1$ because every return preserves $m_1$. The vertical coarse plaques are tame by~\cite[Prop.~4.6 and Def.~4.10]{Entropy}. Theorem~7.2 there gives $h_\xi\le d_\xi$ for every positive class, whence \[
 0=\sum_{[\xi]:\,\xi(H)>0}(d_\xi-h_\xi),\qquad d_\xi-h_\xi\ge0
 \quad\Longrightarrow\quad h_\chi=d_\chi .
\]
In~\cite[Theorem~7.2, pp.~488--489]{Entropy} take $f=a$, the measure $\bar m_1$ and the foliation $W_{\cV}^{[\chi]}$. All its exponents are positive for $a$. The equality $h_\chi=d_\chi$ therefore gives $(\bar m_1)^\zeta_y\sim\vol_{\zeta(y)}$ for $\bar m_1$-almost every $y\in X$.

\par\Needspace{3\baselineskip}\smallskip\noindent\textbf{Identification of the plaque germ.} For $a_\pm=\exp(\pm H_0+\epsilon H_1)$ chosen before Lemma~\ref{coarse}, condition \eqref{eq:coarse-sign} gives
\[
 T_yW^s_{\cV,a_+}(y)\cap T_yW^s_{\cV,a_-}(y)=D_0H_z S_\chi,\qquad y=e_0(z).
\]
Here the stable spaces in the fiber are those computed in the entropy calculation above. By the defining intersection before Lemma~\ref{coarse}, the vertical coarse plaque through $y$ lies locally in both stable plaques for $a_\pm$. It is tangent at $y$ to $D_0H_z S_\chi$ by~\cite[Prop.~4.6 and Def.~4.10]{Entropy}. In the fundamental-domain trivialization fix a Riemannian metric on $M$. In the set $D\subset\widehat X$ where the two derivative norms at zero are bounded, choose $r_0>0$ and $A<\infty$ so that
\[
 D'=\{z\in D:H_z\text{ embeds }B_{r_0}(V),\
                 \sup_{w\in B_{r_0}(V)}\|D_wH_z\|\le A\}
\]
has positive $\widehat m$-measure. These sets exhaust $D$ as $r_0\downarrow0$ and $A\uparrow\infty$ since each chart is analytic near zero. Fix $d\in\{a_+,a_-\}$. For $\widehat m$-almost every $z\in\widehat X$ let $n_1<n_2<\cdots$ be the positive times with $d^{n_i}z\in D'$. For $s\in S_\chi$ in a sufficiently small fixed ball, \eqref{susp} gives
\[
 \operatorname{dist}_{\cV}(d^{n_i}H_z(s),d^{n_i}H_z(0))
 \le A\|\rho(d)^{n_i}s\|\le A' e^{-c n_i}\|s\|
 \qquad(c>0).
\]
The identity $d^{n_i}H_z(s)=H_{d^{n_i}z}(\rho(d)^{n_i}s)$ extends to the stated ball by uniqueness on the intersection of the centered domains, using Lemma~\ref{wall} on $S_\chi$. Distances can be measured on $M$ in the fundamental-domain trivialization of~\cite[\S2.4]{Boundary}. Ergodicity of $d$ gives $n_i/i\to1/\widehat m(D')$, hence $n_{i+1}/n_i\to1$. At a frame over $u\in G_0/\Lambda$, put $\ell(z)=C(1+|\gamma_d(u)|)$, with $e^{\ell(z)}$ bounding the one-step fiber Lipschitz constant. Since $\ell\in L^1(\widehat m)$, Birkhoff's theorem gives $\sum_{j=0}^{n-1}\ell(d^jz)=n\int\ell\,d\widehat m+o(n)$. For $n_i\le n<n_{i+1}$, the distance at time $n$ is bounded by
\[
 A'\|s\|\exp\!\left(-cn_i+\sum_{j=n_i}^{n_{i+1}-1}\ell(d^jz)\right)
 =A'\|s\|\exp(-cn+o(n)).
\]
The bound is uniform in $s$. Shrinking the disk places it in the local stable branch for each of $a_\pm$, by the stable manifold theorem~\cite[Prop.~4.6]{Entropy}.

Write $P_\pm$ for these two embedded local plaques in the fiber at $y$. Choose local defining functions for $P_-$ with independent differentials. On restriction to $T_y P_+$ their common kernel is $T_yP_+\cap T_yP_-=D_0H_z S_\chi$. Select a maximal independent set of these restricted differentials. By the implicit function theorem the corresponding functions cut out in $P_+$ a manifold of dimension $\dim S_\chi$ containing $P_+\cap P_-$ near $y$. Both the coarse plaque and the disk $H_z(S_\chi)$ lie in this manifold, pass through $y$ and have its dimension. Each is open in it near $y$, so their germs agree.
\end{proof}

\begin{proposition}[The law of the displaced center]\label{eval}
With $r$ as in Lemma~\ref{wall} and $S_\chi=\bigoplus_{\lambda=t\chi,\ t>0}V_\lambda$, the distribution of $H_z(s)$ for $z$ sampled from the suspended frame probability $\widehat m$ is
\[
 (e_s)_*\widehat m=(e_0)_*\widehat m=\bar m_1
 \qquad(s\in B_r(S_\chi)).
\]
\end{proposition}
\begin{proof}
Fix a Borel set $D\subset X$ with $\bar m_1(D)=0$ and a partition $\zeta$ subordinate to the plaques of class $\chi$. By Lemma~\ref{coarse},
\[
 0=\int(\bar m_1)^\zeta_y(D)\,d\bar m_1(y)
 \quad\Longrightarrow\quad
 \vol_{\zeta(y)}(D\cap\zeta(y))=0
 \quad(\bar m_1\text{-a.e. }y\in X).
\]
For $\widehat m$-almost every $z\in\widehat X$, let $\delta(z)=2^{-n}r$, where $n$ is the least nonnegative integer with
\[
 H_z\bigl(B_{2^{-n}r}(S_\chi)\bigr)\subset\zeta(e_0(z)).
\]
Lemma~\ref{wall} gives an embedding on this disk, and Lemma~\ref{coarse} and subordination give such an $n$. For fixed $n$, failure of the inclusion is the projection of the Borel set of $(z,s)\in\widehat X\times B_{2^{-n}r}(S_\chi)$ with $H_z(s)\notin\zeta(e_0(z))$, using Borel versions on the regular set. Thus $\delta$ is measurable for the completion of $\widehat m$. This choice depends only on the atom and the chart. Change of variables on the disk gives
$\int_{B_{\delta(z)}(S_\chi)}{\bf1}_{D}(H_z(s))\,ds=0$, with $ds$ Lebesgue measure on $S_\chi$.

For $d=a_+$, condition \eqref{eq:coarse-sign} gives $\|\rho(d)^n|_{S_\chi}\|\to0$ as $n\to\infty$. Choose $j$ with $\widehat m\{z\in\widehat X:\delta(z)\ge1/j\}>0$. Apply the last equality to the countable family $d^nD$, $n\ge0$. For $\widehat m$-almost every $z\in\widehat X$ choose a return $n$ with $\delta(d^nz)\ge1/j$ and $\rho(d)^n B_r(S_\chi)\subset B_{1/j}(S_\chi)$. Such returns exist by ergodicity of $d$. On $B_r(S_\chi)$, formula~\eqref{susp} continues along each segment $[0,s]$ by Lemma~\ref{wall}. Substitution $v=\rho(d)^ns$ gives
\[
 \int_{B_r(S_\chi)}{\bf1}_D(H_z(s))\,ds
 =|\det(\rho(d)^n|_{S_\chi})|^{-1}
 \int_{\rho(d)^nB_r(S_\chi)}{\bf1}_{d^nD}(H_{d^nz}(v))\,dv=0 .
\]
For $B=B_r(S_\chi)$ let $|B|$ denote its Lebesgue volume. Fubini gives
\begin{equation}\label{av}
 \frac1{|B|}\int_B(e_s)_*\widehat m\,ds\ll\bar m_1 .
\end{equation}
Each probability in this integral is $b$-invariant by Lemma~\ref{wall}. The density of the average relative to $\bar m_1$ is $b$-invariant. It is constant by ergodicity and equals one by normalization. Extremality of $\bar m_1$ gives $(e_s)_*\widehat m=\bar m_1$ for Lebesgue-almost every $s\in B$. For any bounded continuous $\varphi$ on $X$, dominated convergence gives continuity of $s\mapsto\int\varphi(H_z(s))\,d\widehat m(z)$. Its equality with $\int\varphi\,d\bar m_1$ on a set of full Lebesgue measure therefore holds on all of $B$.
\end{proof}

\section{The two-point joining}

\begin{lemma}[Change of coordinates at a fixed displacement]\label{joining}
Let $b=\exp H_0$ be ergodic on $(\widehat X,\widehat m)$, $C_b=\Fix(\rho(b))$, and use the common ball of Lemma~\ref{wall}. If $s\in B_r(C_b)$ satisfies $(e_s)_*\widehat m=\bar m_1$, there is an invertible measure-preserving transformation $T_u$ of $\widehat X$ modulo null sets and one analytic local diffeomorphism $u$ satisfying
\[
 u(0)=s,\qquad u\rho(b)=\rho(b)u,\qquad H_{T_u z}=H_z\circ u
 \quad(\widehat m\text{-almost every }z\in\widehat X).
\]
\end{lemma}
\begin{proof}
Fix $s\in B_r(C_b)$. For $v=0,s$, the evaluation of \eqref{susp} is $e_v:\widehat X\to X$, $e_v(z)=H_z(v)$, defined on a set of full $\widehat m$-measure by Lemma~\ref{wall}. For $\bar m_1$-almost every $y\in X$ let $\widehat m_y^v$ be the conditional probability on $e_v^{-1}(y)$. The hypothesis $(e_s)_*\widehat m=(e_0)_*\widehat m=\bar m_1$ gives the same measure on $X$ in both disintegrations. Thus the probability
\begin{equation}\label{join}
 J_s=\int_X \widehat m_y^s\otimes\widehat m_y^0\,d\bar m_1(y)
 \quad\text{satisfies}\quad J_s\{(z,z')\in\widehat X^2:H_z(s)=H_{z'}(0)\}=1 .
\end{equation}
Both evaluations commute with $b$, so $(b,b)_*J_s=J_s$. Write $J_s=\int\xi\,d\tau_s(\xi)$ for its ergodic decomposition under $(b,b)$. Let $\operatorname{pr}_1,\operatorname{pr}_2:\widehat X^2\to\widehat X$ be the coordinate projections. Both marginals of $J_s$ are the $b$-ergodic probability $\widehat m$. Extremality applied to each marginal, and \eqref{join} applied to the incidence set, give
\[
 (\operatorname{pr}_1)_*\xi=(\operatorname{pr}_2)_*\xi=\widehat m,\qquad
 \xi\{(z,z')\in\widehat X^2:H_z(s)=H_{z'}(0)\}=1
\]
for $\tau_s$-almost every $\xi$. Fix one such component.

\par\Needspace{3\baselineskip}\smallskip\noindent\textbf{The change of coordinates.}
For $\xi$-almost every $(z,z')\in\widehat X^2$, we have $H_{z'}(0)=H_z(s)$. Use the inverse branch of $H_z$ taking this point to $s$. Then $u_{z,z'}=H_z^{-1}\circ H_{z'}:(V,0)\to(V,s)$ is an analytic local diffeomorphism. Formula~\eqref{susp} and $\rho(b)s=s$ give, for $A_j(z,z')=D^j_0u_{z,z'}/j!$,
\[
 A_j(bz,bz')=\rho(b)A_j(z,z')(\rho(b)^{-1})^{\otimes j},\qquad j\ge1 .
\]
Here $V_\lambda=\{v\in V:\rho(\exp H)v=e^{\lambda(H)}v\ (H\in\mathfrak a)\}$ is the weight space used in Lemma~\ref{coarse}. Project the output onto $V_{\lambda_0}$ and restrict input $i$ to $V_{\lambda_i}$ for $1\le i\le j$. For the resulting component $A$ the multiplier is $a=\exp((\lambda_0-\sum_{i=1}^j\lambda_i)(H_0))$. If $a\ne1$, put $q=\max(a,a^{-1})>1$. Invariance of $\xi$ under $(b,b)^{\pm1}$ gives, for every $R>0$,
\[
 \xi\{(z,z')\in\widehat X^2:\|A(z,z')\|>R\}
 =\xi\{(z,z')\in\widehat X^2:\|A(z,z')\|>q^nR\}\longrightarrow0 .
\]
Thus $A=0$ $\xi$-almost everywhere when $a\ne1$. For $a=1$, ergodicity of $\xi$ makes $A$ constant. Intersect these sets of full $\xi$-measure over the finitely many components at each $j$ and the countably many $j$. At one pair in this intersection the common Taylor series represents an analytic local diffeomorphism $u$. At every pair in the intersection, analyticity gives $H_{z'}=H_z\circ u$ as germs at zero. Coefficient invariance gives $u\rho(b)=\rho(b)u$.

\par\Needspace{3\baselineskip}\smallskip\noindent\textbf{The map of frames.}
For fixed $z\in\widehat X$, the center $H_z(s)$ and derivative $D_sH_z\,D_0u$ determine the terminal frame $z'$ by \eqref{frames} and \eqref{susp}. With the terminal frame fixed, $H_z=H_{z'}\circ u^{-1}$ near $s$. Any two candidates are analytic near $[0,s]$ by Lemma~\ref{wall}. The component of their common domain that contains this segment also contains a neighborhood of $0$. Hence their germs at $0$ agree, and they give the same frame. The two disintegrations of $\xi$ over its marginals are therefore point~masses. They give inverse measurable maps, one of which satisfies
\[
 \xi=(z\mapsto(z,T_u z))_*\widehat m,\qquad
 (T_u)_*\widehat m=\widehat m,\qquad H_{T_u z}=H_z\circ u
 \quad(\widehat m\text{-a.e. }z\in\widehat X).
\]
If a second component in this full $\tau_s$-measure set gives $\tilde u$ over $s$, then $H_{T_uz}(0)=H_z(s)=H_{T_{\tilde u}z}(0)$ for $\widehat m$-almost every $z\in\widehat X$. By \eqref{frames} and centralization of $\rho(G_0)$ by $K$, their charts differ on the right by $h\in K$. Cancel $H_z$ near $s$ to obtain $\tilde u=u\circ h$. In particular $h=(D_0u)^{-1}D_0\tilde u$ is independent of $z$.
\end{proof}

\begin{lemma}[One-parameter changes of coordinates]\label{lie}
Keep $b$ and $C_b$ of Lemma~\ref{joining}. Let $S\subset C_b$ be a linear subspace and assume $(e_s)_*\widehat m=\bar m_1$ for all sufficiently small $s\in S$, where $e_s(z)=H_z(s)$ by Lemma~\ref{wall}. For every $v\in S$ there is an analytic one-parameter family of local changes $u_t$ as in Lemma~\ref{joining}, with
\[
 u_0=\id,\qquad \left.\frac{d}{dt}\right|_0u_t(0)=v .
\]
\end{lemma}
\begin{proof}
Fix $b$ and put $C_b=\{w\in V:\rho(b)w=w\}$. Consider analytic local diffeomorphisms $u$ from a neighborhood of $0$ in $V$ to $V$, with $u(0)$ near $0$ in $C_b$, satisfying $u\rho(b)=\rho(b)u$ and admitting an invertible $\widehat m$-preserving $T_u$ with $H_{T_uz}=H_z\circ u$ for $\widehat m$-almost every $z\in\widehat X$, as in Lemma~\ref{joining}. The equation $u\rho(b)=\rho(b)u$ gives $u(C_b)\subset C_b$ wherever both sides are defined near $0$.

\par\Needspace{3\baselineskip}\smallskip\noindent\textbf{Common domains.}
Choose a measurable set $D\subset\widehat X$ with $\widehat m(D)>1/2$ on which the charts $H_z$ and their inverse branches have bounded holomorphic extensions on balls of a common radius. To obtain $D$, use a measurable section of $G_0\to G_0/\Lambda$ and finitely many relatively compact analytic charts of the compact fiber $M$. At each center some positive rational radius and some integer bound work for both representatives. Smaller coordinate neighborhoods covering $M$ allow the two representatives at sufficiently close centers to be expressed in a common chart with bounded transition maps. Taking a sufficiently small radius and a sufficiently large bound gives the required measure.

For each $u$, measure preservation gives
\[
 \widehat m(D\cap T_u^{-1}D)\ge 2\widehat m(D)-1>0 .
\]
Choose $z$ there where $H_{T_uz}=H_z\circ u$. The two centers in the same fiber are $H_z(0)$ and $H_z(u(0))$. Their distance is at most $C\|u(0)\|$ by the uniform derivative bound. For $\|u(0)\|$ sufficiently small, each center lies inside the other's inverse domain with a common positive distance to its boundary. Shrinking the source ball using the derivative bound gives uniform holomorphic bounds on a fixed ball for
\[
       u=H_z^{-1}\circ H_{T_u z},\qquad u^{-1}=H_{T_u z}^{-1}\circ H_z .
\]
Writing $s=u(0)$, the second formula gives $u^{-1}(s)=0$. Both $H_z(w)$ and $H_{T_uz}^{-1}(H_z(w))$ are defined for $w$ on a common ball containing $[0,s]$. If $C_1$ bounds $\|Du^{-1}\|$ there, integration along the segment gives $\|u^{-1}(0)\|\le C_1\|s\|$. Write $B_R^{\mathbb C}=\{w\in V\otimes_{\mathbb R}\mathbb C:\|w\|<R\}$ for $R>0$, using the Hermitian extension of the inner product on $V$. By the common holomorphic bounds and Cauchy's estimate, choose $R$ and $A\ge1$ so that $Du$ and $Du^{-1}$ have norm at most $A$ on $B_{2R}^{\mathbb C}$ for every such change. Take $2R$ smaller than the radius in Lemma~\ref{wall}, put $\delta=R/(4A)$, and restrict the displacements so that $\|u(0)\|,\|u^{-1}(0)\|<\delta$. Then
\[
 \|u^{\pm1}(w)\|\le\delta+A\|w\|<R/2
 \qquad(w\in B_\delta^{\mathbb C}).
\]
Thus every composition of two such maps or inverses is defined on $B_\delta^{\mathbb C}$.

\par\Needspace{3\baselineskip}\smallskip\noindent\textbf{Composition and limits.}
For two changes $u,v$ with these displacements, $[0,v(0)]$ and its image under $u$ lie in the common ball of $C_b=\Fix(\rho(b))$ from Lemma~\ref{wall}. Continue $H_{T_u z}=H_z\circ u$ along $[0,v(0)]$ in $C_b$, and compose at its endpoint with $v$. This gives
\[
 H_{T_vT_uz}=H_z\circ u\circ v,\qquad T_{u\circ v}=T_vT_u
 \quad\bmod\widehat m.
\]
The inverse commutes with $\rho(b)$ near $u(0)$. Restriction to the fixed space and continuation along $[u(0),0]$ give $u^{-1}(0)\in C_b$. Continuing the frame identity along $[0,u^{-1}(0)]$ now gives $T_{u^{-1}}=T_u^{-1}$.

Use locally uniform convergence on smaller complex balls for maps and inverses. For any sequence of these changes with $u_n(0)$ in a sufficiently small ball of $C_b$, Montel's theorem and the common bounds give a subsequence converging to $(u,v)$ for the maps and inverses. The identity $u_n^{-1}\circ u_n=\id$ holds near $0$, and $u_n\circ u_n^{-1}=\id$ holds near $u_n(0)$. Both extend to $B_\delta^{\mathbb C}$, a connected ball containing both points on which the compositions are defined. With $\|f\|_D=\sup_{w\in D}\|f(w)\|$, we have
\[
 \|u_n\circ u_n^{-1}-u\circ v\|_{B_\delta^{\mathbb C}}
 \le\|u_n-u\|_{B_R^{\mathbb C}}+A\|u_n^{-1}-v\|_{B_\delta^{\mathbb C}}\longrightarrow0.
\]
The same estimate in the opposite order gives $u\circ v=v\circ u=\id$ on $B_\delta^{\mathbb C}$. Let $J(z)=j^\infty_0H_z$, the sequence of all Taylor coefficients, including the center. In local fiber coordinates this sequence belongs to a countable product of finite-dimensional jet spaces. The map $J$ is Borel and injective by \eqref{frames}, and its inverse on its image is Borel~\cite[Corollary~15.2]{Kechris}. Since $u_n(0)\to u(0)$ in $C_b$, Lemma~\ref{wall} gives for each order convergence of the jets of $H_z\circ u_n$ at zero for $\widehat m$-almost every $z\in\widehat X$. For $z=[g,p]$, the order-$j$ jet includes its center in the fiber of $X$ over $g\Lambda$. Its lower-order jets are its truncations. Testing bounded continuous functions at each finite order therefore gives equality on the cylinder sets generating the Borel sets of the jet sequence. It follows that
\[
 (z\mapsto j^\infty_0(H_z\circ u))_*\widehat m=J_*\widehat m .
\]
Consequently $T_u z=J^{-1}(j^\infty_0(H_z\circ u))$ exists almost everywhere and preserves $\widehat m$. Apply the same argument to $u^{-1}$ and use the composition identity to obtain its inverse. Uniform bounds, uniqueness of analytic continuation and Montel's theorem give a compact closed neighborhood of the identity in this topology by imposing closed small bounds on the displacements and on the sup norms of $u-\id$ and $u^{-1}-\id$ on a smaller ball. Cauchy's estimates give continuity of composition and inversion after a further shrinking. Uniform convergence on any such smaller ball gives the same topology on this bounded family, by Montel's theorem and uniqueness.

\par\Needspace{3\baselineskip}\smallskip\noindent\textbf{One-parameter subgroups.}
Choose an identity neighborhood on which $\|Du-I\|<1/2$ on a common ball. Let $H$ be the compact closure of a subgroup contained in a smaller neighborhood. With $m_H$ its Haar probability, put $A(w)=\int_H h(w)\,dm_H(h)$. On a smaller ball $\|DA-I\|\le1/2$, so $A$ is injective. For $h\in H$, right invariance gives $A\circ h=A$ wherever both sides are defined, hence $h=\id$ as a germ. Thus the locally compact local group has no small subgroups. By~\cite[\S8]{Goldbring} an identity neighborhood is a local Lie group.

For a continuous one-parameter subgroup $u_t$ and $\epsilon>0$ small put $B_\epsilon(w)=\epsilon^{-1}\int_0^\epsilon u_t(w)\,dt$. The common holomorphic bounds give analyticity in $w$, and $D_0B_\epsilon\to I$. The group law gives
\[
 B_\epsilon(u_s(w))=\epsilon^{-1}\int_s^{s+\epsilon}u_t(w)\,dt,\qquad
 \left.\frac{d}{ds}\right|_{s=0}u_s(w)
   =(D_wB_\epsilon)^{-1}\frac{u_\epsilon(w)-w}{\epsilon}.
\]
Here $D_wB_\epsilon:V\to V$ is the derivative in $w$ at fixed $\epsilon$. The right side is an analytic vector field. The group law and uniqueness of its local flow give joint analyticity in $(s,w)$. Products of these flows in exponential coordinates give an analytic local group action.

If $u(0)=0$, the two frames have the same center, so cancellation gives $u=h$ for some $h\in K$. Conversely $z\mapsto zh$ realizes each $h\in K$ by \eqref{frames}. Thus the stabilizer of $0$ is $K$. For every sufficiently small $s\in S$, Lemma~\ref{joining} gives a change $u$ with $u(0)=s$ and $H_{T_uz}=H_z\circ u$ for $\widehat m$-almost every $z\in\widehat X$. For every neighborhood $\mathcal U$ of the identity there is $\delta>0$ such that $B_\delta(S)\subset\{u(0):u\in\mathcal U\}$. Otherwise choose omitted $s_n\to0$ and changes $u_n(0)=s_n$. Compactness for bounded small displacements gives a subsequence $u_n\to h\in K$. Then $u_nh^{-1}\to\id$ and $u_nh^{-1}(0)=s_n$, a contradiction. Left translation in the local group and invertibility of the changes give constant rank for $u\mapsto u(0)$. Its local image contains a neighborhood of $0$ in $S$, so its tangent contains $S$. For $v\in S$, choose a Lie algebra element mapped to $v$ by its differential at the identity. Its one-parameter subgroup has the required velocity.
\end{proof}

\section{Analytic local sections}\label{sec:analytic-sections}

\begin{proof}[Proof of Lemma~\ref{sections}]
Recall from \eqref{susp} that $H_z=\iota_g\psi_p\rho(g)^{-1}$ for $z=[g,p]\in\widehat X$, with $\iota_g(x)=[g,x]\in X$. Write $V=U\oplus F$, where $F=V^{G_0}$ and $U$ is the sum of the nontrivial isotypic summands. The $G_0$-translates of the positive coarse spaces $S_\chi=\bigoplus_{t>0}V_{t\chi}$, over all nonzero weights $\chi$, span $U$. Indeed their span is invariant and lies in $U$. On the quotient of $U$ by that span the split torus acts trivially. The kernel of the representation is normal and contains noncentral elements in each simple factor, so the whole quotient is trivial as a representation. Complete reducibility of $U$ forces the quotient to be zero.

By Proposition~\ref{eval} and Lemma~\ref{lie}, the one-parameter changes have velocities at zero spanning each $S_\chi$. If $u_t$ is one of these changes and $d\in G_0$, formula~\eqref{susp} gives the change $\rho(d)u_t\rho(d)^{-1}$ with transformation $dT_{u_t}d^{-1}$. Its center path lies in $\rho(d)C_b=\Fix(\rho(dbd^{-1}))$. Choose families $u^1,\ldots,u^\ell$ whose velocities $v_i=\partial_t|_0u^i_t(0)$ form a basis of $U$. For $t=(t_1,\ldots,t_\ell)$ near $0$, the changes and the frame transformations are composed in opposite orders,
\[
 u(t)=u^1_{t_1}\circ\cdots\circ u^\ell_{t_\ell},\qquad
 T^{(i)}(t)=T_{u^i_{t_i}}\circ\cdots\circ T_{u^1_{t_1}},
 \qquad T^{(0)}(t)=\id .
\]
Here $u(t)$ is a local map in $V$, whereas $T^{(i)}(t)$ acts on $\widehat X=(G_0\times P)/\Lambda$. For each fixed $t$, every $T^{(i)}(t)$ preserves $\widehat m$, the suspension of the frame probability $\eta$ in \eqref{susp}. To obtain joint measurability in $t$ and $z$ for one factor, use the Borel inverse of $z\mapsto j^\infty_0H_z$ in Lemma~\ref{lie} at $j^\infty_0(H_z\circ u^i_{t_i})$ on its domain. Compose these versions and assign arbitrary values on their exceptional sets.

\begin{claim}\label{cl:fixed-evaluation}
For some $r>0$, $(e_f)_*\widehat m=\bar m_1$ for every $f\in B_r(F)$, where $F=V^{G_0}=\{v\in V:\rho(g)v=v\ (g\in G_0)\}$ and $e_f([g,p])=[g,\psi_p(f)]$ by \eqref{susp}.
\end{claim}
\begin{claimproof}
Choose $H\in\mathfrak a$ with nonzero projection to every simple factor and let $b=\exp H$. Since $F\subset\Fix(\rho(b))$, Lemma~\ref{wall} gives $r>0$ such that $H_z$ embeds a neighborhood in $V$ of $\overline B_{2r}(F)$ for $\widehat m$-almost every $z\in\widehat X$. If $U=0$, then $F=V$ and the empty composition $u(t)$ is the identity. If $F=0$, the claim is $(e_0)_*\widehat m=\bar m_1$, already proved after~\eqref{susp}. Use $\delta_0$ for Lebesgue probability on a zero-dimensional parameter space. Shrink $r$ so the finitely many families $u^i_t$ are defined there for all sufficiently small $t$. The map $(t,f)\mapsto u(t)(f)$ has derivative at $(0,0)$
\[
  (\dot t,\dot f)\longmapsto\sum_{i=1}^\ell\dot t_i v_i+\dot f.
\]
Since $V=U\oplus F$, it is a diffeomorphism on a product of neighborhoods. Fix $B_r(F)$ inside this product for the rest of the proof.

For $\delta>0$ let
\[
 \begin{aligned}
 N_\delta&=\{w\in V:\operatorname{dist}(w,\overline B_r(F))<\delta\},\\
 D_\delta&=\{z\in\widehat X:H_z\text{ has an embedding representative on }N_\delta\}.
 \end{aligned}
\]
Given $\epsilon>0$, choose $\delta>0$ with $\widehat m(D_\delta)>1-\epsilon$. This is possible by compactness of $\overline B_r(F)$ and measurability of the charts. Choose a time box $I$ about $0$ so small that all partial compositions of consecutive factors move points of $\overline B_r(F)$ by less than $\delta/2$. For $t\in I$ put
\[
 D(t)=\bigcap_{i=0}^{\ell}(T^{(i)}(t))^{-1}D_\delta,\qquad
 \widehat m(D(t)^c)\le(\ell+1)\epsilon .
\]
For each fixed $t$, the one-step identities of Lemma~\ref{joining} hold at every successive frame for $\widehat m$-almost every $z\in\widehat X$. For such $z\in D(t)$ continue these identities from $0$ along the segments in $\overline B_r(F)$ and their images under the partial compositions. Indeed we may shrink $\delta$ so the individual families are defined on the $\delta$-tube about $\overline B_r(F)$, then shrink $I$ so all consecutive compositions differ from the identity by less than $\delta/4$ on the closed $\delta/2$-tube. The latter tube is connected and lies inside each of the required composition domains. Successive substitution gives
\begin{equation}\label{words}
 H_{T^{(\ell)}(t)z}(f)=H_z(u(t)(f))\qquad(f\in B_r(F)).
\end{equation}
Let $A\subset X$ be Borel with $\bar m(A)=0$. For Haar-almost every base coset, its section in $M$ has $m$-measure zero. Fix such a coset $g\Lambda$ and an initial $z=[g,p]$. The map $(t,f)\mapsto H_z(u(t)(f))$ is a local diffeomorphism into that fiber wherever $H_z$ embeds at $u(t)(f)$. Change of variables therefore gives
\[
 \mathcal L^{\dim V}\{(t,f)\in I\times B_r(F):z\in D(t),\ H_z(u(t)(f))\in A\}=0
\]
for $\widehat m$-almost every $z\in\widehat X$. Normalize Lebesgue measure on $I$ and $B_r(F)$ to probabilities $dt$ and $df$. Fubini, \eqref{words}, and $(T^{(\ell)}(t))_*\widehat m=\widehat m$ give
\[
 \begin{aligned}
 \int_{B_r(F)} (e_f)_*\widehat m(A)\,df
 &=\int_{I\times B_r(F)}\int {\bf1}_A(H_{T^{(\ell)}(t)z}(f))\,d\widehat m(z)\,dt\,df\\
 &\le \int_I\widehat m(D(t)^c)\,dt\le(\ell+1)\epsilon .
 \end{aligned}
\]
Only $I$ and $\delta$ depend on $\epsilon$. Thus the average in $f$ is absolutely continuous with respect to $\bar m$.

Fix $f\in B_r(F)$. Its defining property is $\rho(g)f=f$ for every $g\in G_0$. By \eqref{susp} and continuation along $[0,f]$, $e_f(gz)=g e_f(z)$ for $\widehat m$-almost every $z\in\widehat X$ at each fixed $g\in G_0$. Its law is therefore $G_0$-invariant and ergodic. Since $\bar m=k^{-1}\sum_j\bar m_j$, an invariant probability absolutely continuous with respect to $\bar m$ is a convex combination of these $k$ components because its density is invariant and hence constant on each. Uniqueness of ergodic decomposition applied to this average now gives $(e_f)_*\widehat m\in\{\bar m_1,\ldots,\bar m_k\}$ for Lebesgue-almost every $f\in B_r(F)$. The law is weakly continuous in $f$ by dominated convergence. A set of full Lebesgue measure is dense in $B_r(F)$, and a continuous image of this connected ball in a finite set is a point. Evaluation at $0$ gives $\bar m_1$.
\end{claimproof}

Apply Lemmas~\ref{joining} and~\ref{lie} to the fixed space in Claim~\ref{cl:fixed-evaluation}. Append families whose velocities form a basis of $F$. Use the same order as above for the resulting $d=\dim V$ families, and write $u(t)$ and $T^{(d)}(t)$ for their compositions. The derivative at $0$ of $t\mapsto u(t)(0)$ sends the standard basis to the chosen basis of $V$.

\begin{claim}\label{cl:section-word}
For $\widehat m$-almost every $z\in\widehat X$ there is a box $I_z$ about $0$ such that, for Lebesgue-almost every $t\in I_z$,
\[
 H_{T^{(d)}(t)z}=H_z\circ u(t)\quad\text{as germs at }0,\qquad
 H_z(u(t)(0))\in(G_0\times X_1)/\Lambda .
\]
\end{claim}
\begin{claimproof}
Fix $R>0$ and take initial frames $z$ for which $H_z$ embeds $B_R(V)$. Choose a time box so small that every prefix of $u(t)$ and every image, under the preceding prefix, of an individual center path $s\mapsto u^i_s(0)$ stays in $B_R(V)$ for $s$ between $0$ and $t_i$. For fixed $t$, the inverse images of the exceptional sets under the preceding frame transformations have $\widehat m$-measure zero. Fix an initial $z$ outside this null set for which $H_z$ embeds $B_R(V)$. Suppose the identity has been proved for the first $i-1$ factors. Put $z_{i-1}=T^{(i-1)}(t)z\in\widehat X$ and $v=u^1_{t_1}\circ\cdots\circ u^{i-1}_{t_{i-1}}$ (the identity if $i=1$). The induction hypothesis gives $H_{z_{i-1}}=H_z\circ v$ near $0$. For $a$ between $0$ and $t_i$, the path $u^i_a(0)$ lies in the fixed space of the Cartan element for this family. Lemma~\ref{wall} gives $H_{z_{i-1}}$ on a neighborhood of this path; the choice of the box puts $v(u^i_a(0))$ in $B_R(V)$. Thus the common domain of $H_{z_{i-1}}$ and $H_z\circ v$ contains the path. Analytic continuation gives equality near its endpoint $u^i_{t_i}(0)$. Composing at that point gives
\[
 H_{T^{(i)}(t)z}
 = H_{z_{i-1}}\circ u^i_{t_i}
 = H_z\circ v\circ u^i_{t_i}
 \quad\text{as germs at }0 .
\]
This proves the first identity for $i=d$.
The second holds for this fixed $t$ because $T^{(d)}(t)$ preserves $\widehat m$ and $m_1(X_1)=1$, where $m_1=k\,m|_{X_1}$ by \eqref{components}. Apply Fubini on each box and take the union over $R=1/n$, $n\ge1$.
\end{claimproof}

Fix $z=[g,p]$ as in Claim~\ref{cl:section-word}. Recall $\iota_g(x)=[g,x]$ from \eqref{susp}. On a smaller time box the center map $t\mapsto\iota_g^{-1}H_z(u(t)(0))$ is an analytic diffeomorphism onto an open neighborhood $U_0$ of $c(p)$. Denote its inverse by $x\mapsto t(x)$. The section on $U_0$ is
\[
 s_x(w)=\iota_g^{-1}H_z\!\left(u(t(x))(\rho(g)w)\right).
\]
It is analytic in $(x,w)$ and $s_x(0)=x$. If $T^{(d)}(t)z=[g,q]$ and the word identity holds, \eqref{susp} gives $H_{[g,q]}(\rho(g)w)=\iota_g(\psi_q(w))$, hence $s_x=\psi_q$. Change of variables from $t$ to $x$ gives this equality for $m$-almost every $x\in U_0$ and gives $m(U_0\setminus X_1)=0$.

Let $O_1=\bigcup_{\gamma\in\Lambda}\gamma U_0$. By \eqref{components}, $m_1(U_0)=k\,m(U_0)>0$. Ergodicity gives $m_1(O_1)=1$ and $m_1|_{O_1}=k\,m|_{O_1}$. On $\gamma U_0$ take the section $w\mapsto\gamma s_{\gamma^{-1}x}(\rho(\gamma)^{-1}w)$. For two sections $s_x,t_x$ on an overlap, put $h(x)=(D_0s_x)^{-1}D_0t_x$. At $m$-almost every $x$, $h(x)\in K$ and $t_x=s_x\circ h(x)$ by \eqref{frames}. Since $K$ is closed, the coefficients are continuous in $x$, and $m$ is positive on every nonempty open set, both identities hold throughout the overlap. These analytic transition functions define the bundle, and the conditional law on each $K$-fiber is $m_K$ by Proposition~\ref{action-coordinates}.
\end{proof}

\manuscriptpart{Homogeneous realization}\label{assembly}
\noindent\textbf{Outline of Part III.}
Lemma~\ref{finite} lifts every element stabilizing $m_1$ to the analytic frame bundle. The first derivative of a chart gives its frame, and the second derivative gives the connection term in Lemma~\ref{coframe}. Ergodicity makes the brackets of the dual fields constant in this coframe. The fields initially have local flows. Theorem~\ref{H-thm:coframe-integration} uses elements of the split torus to prove the group relations for their measure-preserving extensions. Lemma~\ref{H-lem:coframe-measured-realization} identifies a local flow through a frame with right multiplication on its orbit. The proof of Theorem~\ref{main} divides these orbits by $K$ and puts back the finitely many translates of $m_1$.

\section{The stabilizer of the volume component}
Let $\Gamma_1=\{t\in\Gamma:t_*m_1=m_1\}$, so $\Lambda\lhd\Gamma_1$ and $[\Gamma:\Gamma_1]=k$ by \eqref{components}. Set $O_1'=\bigcap_{t\in\Gamma_1}tO_1$. This is a finite intersection because $\gamma O_1=O_1$ for $\gamma\in\Lambda$ and $\Lambda\lhd\Gamma_1$. By Lemma~\ref{sections} this open set is $\Gamma_1$-invariant and $m_1(O_1')=1$. Write $c:P_1\to O'_1$ for the restricted frame bundle. In every section of Lemma~\ref{sections}, $d\eta=dm_1\,dm_K$. Thus $\eta$ has full support on $P_1$.
\begin{lemma}[Lifting the stabilizer of the component]\label{finite}
In a single analytic source coordinate commuting with $\rho(G_0)$ and $K$, each $t\in\Gamma_1$ has an analytic measure-preserving lift $T_t$ to this bundle and a linear map $A_t$ normalizing $K$ such that
\[
 t\psi_p=\psi_{T_tp}A_t,\qquad T_tT_s(p)=T_{ts}(p)k_{t,s}\quad(k_{t,s}\in K).
\]
The source changes $A_t$ belong to the product of $\rho(G_0)$ and one compact linear group normalizing $K$ and centralizing $\rho(G_0)$.
\end{lemma}

\begin{proof}
\textbf{Comparison with a translate.}
Use $\gamma\psi_p=\psi_{\gamma p}\rho(\gamma)$ from \eqref{frames} for $\gamma\in\Lambda$.
For $t\in\Gamma_1$ choose a lift $(g_t,c_t)$ to the covering
lattice in Proposition~\ref{prop:group-reduction} and set $R_t=\rho(g_t)$.
Then
\[
 R_t\rho(\gamma)R_t^{-1}=\rho(t\gamma t^{-1})
 \quad(\gamma\in\Lambda).
\]
Fix $t\in\Gamma_1$ and disintegrate $\eta=\int\eta_x\,dm_1(x)$ over $c$. For the probability $J_t=\int\eta_x\otimes\eta_{tx}\,dm_1(x)$ on $P_1\times P_1$, we have $J_t\{(p,q)\in P_1^2:c(q)=tc(p)\}=1$. Thus equivariance and $(\gamma)_*m_1=m_1$ give $(\gamma,t\gamma t^{-1})_*J_t=J_t$. Write $J_t=\int\xi\,d\tau_t(\xi)$ for the ergodic decomposition under $\gamma\cdot(p,q)=(\gamma p,t\gamma t^{-1}q)$. Since both marginals are the ergodic probability $\eta$, extremality gives $(\operatorname{pr}_1)_*\xi=(\operatorname{pr}_2)_*\xi=\eta$ for $\tau_t$-almost every $\xi$, where $\operatorname{pr}_i:P_1^2\to P_1$ are the coordinate projections. The identity for $J_t$ above gives $\xi\{(p,q)\in P_1^2:c(q)=tc(p)\}=1$ on the same full-measure set. Fix such a component. For $\xi$-almost every $(p,q)\in P_1^2$ define
\[
 C_{p,q}=\psi_q^{-1}t\psi_pR_t^{-1}:(V,0)\longrightarrow(V,0).
\]
Put $\delta=t\gamma t^{-1}$. Substitution of $\psi_{\delta q}=\delta\psi_q\rho(\delta)^{-1}$, $\psi_{\gamma p}=\gamma\psi_p\rho(\gamma)^{-1}$ and $R_t\rho(\gamma)R_t^{-1}=\rho(\delta)$ gives
\[
 C_{\gamma p,\delta q}
 =\rho(\delta)\psi_q^{-1}\delta^{-1}t\gamma\psi_p\rho(\gamma)^{-1}R_t^{-1}
 =\rho(\delta)C_{p,q}\rho(\delta)^{-1}.
\]
At degree $n$ this is the action on $\Hom(\Sym^n V,V)$ given by
\[
 A\longmapsto\rho(t\gamma t^{-1})A
       (\rho(t\gamma t^{-1})^{-1})^{\otimes n}.
\]
For each $n\ge1$, put $[C_{p,q}]_n=D^n_0 C_{p,q}/n!$. The coefficient map is equivariant for this action. Since $t\Lambda t^{-1}=\Lambda$ and $\xi$ is ergodic, Lemma~\ref{H-lem:fixed-law} applies. It gives a $\rho(G_0)$-fixed coefficient $c_n$ with $\xi\{(p,q)\in P_1^2:[C_{p,q}]_n=c_n\}=1$. Hence
\[
 \xi\{(p,q)\in P_1^2:[C_{p,q}]_n=c_n\text{ for every }n\ge1\}=1.
\]
Choose $(p,q)\in P_1^2$ in this set. Analyticity gives $\sum_{n\ge1}\|c_n\|r^n<\infty$ for some $r>0$. Thus
\[
 C_t(v)=\sum_{n\ge1}c_n(v^{\otimes n})\quad(\|v\|<r),\qquad B_t=C_tR_t.
\]
For $\xi$-almost every $(p,q)\in P_1^2$, $\psi_q=t\psi_pB_t^{-1}$ and $\psi_p=t^{-1}\psi_qB_t$ as germs. Proposition~\ref{action-coordinates} determines each frame from its center and first derivative. Both directions are Borel by the Borel inverse on the image of the injective first-frame map. Thus $\xi=(p\mapsto(p,T_tp))_*\eta$, $(T_t)_*\eta=\eta$, and for $\eta$-almost every $p\in P_1$,
\begin{equation}\label{H-eq:finite-lift}
 t\psi_p=\psi_{T_tp}B_t,\qquad
 T_t(\gamma p)=(t\gamma t^{-1})T_t p.
\end{equation}
For $\gamma\in\Lambda$ choose $T_\gamma p=\gamma p$ and
$B_\gamma=\rho(\gamma)$.

\textbf{The compact group of source germs.}
Since $cT_t=t c$, one has $T_t(c^{-1}(x))=c^{-1}(tx)$ modulo the conditional null sets over $m_1$-almost every $x\in M$.
For fixed $p\in P_1$ and $h\in K$, the source changes for $p$ and $ph$ give
$\psi_{T_t(ph)}=\psi_{T_tp}\circ B_t h B_t^{-1}$
whenever \eqref{H-eq:finite-lift} holds at both frames. The two charts on the right and left have the same center $tc(p)$. Thus
\[
 B_t K B_t^{-1}=K.
\]
Indeed, for this fixed $t$, substitution gives $B_tkB_t^{-1}\in K$ for $m_K$-almost every $k\in K$. Closedness of $K$ and continuity extend the inclusion to all $k$. Its image contains $K^0$ because it has dimension $\dim K$, and meets every component because it has $|\pi_0(K)|$ components. Hence it equals $K$.
Write $T_tT_s(p)=T_{ts}(p)k(p)$ on the common frame coset. For each $\gamma\in\Lambda$, \eqref{H-eq:finite-lift} gives $k(\gamma p)=k(p)$ for $\eta$-almost every $p\in P_1$. Ergodicity of $\Lambda\curvearrowright(P_1,\eta)$ gives $k_{t,s}\in K$ with $\eta\{p\in P_1:k(p)=k_{t,s}\}=1$. Substitution in \eqref{H-eq:finite-lift} gives
\[
 \psi_{T_{ts}p}B_{ts}=ts\psi_p
 =\psi_{T_tT_sp}B_tB_s=\psi_{T_{ts}p}k_{t,s}B_tB_s .
\]
Cancel the coordinate germ to obtain $B_tB_s=k_{t,s}^{-1}B_{ts}\in K B_{ts}$.
Choose a finite set $\mathcal T\subset\Gamma_1$ with $\Gamma_1=\bigsqcup_{t\in\mathcal T}\Lambda t$,
including $e$, and choose the other lifts by
$B_{\gamma t}=\rho(\gamma)B_t$.
If $ts=\gamma r$ for representatives $t,s,r$, their chosen covering
lifts satisfy
\[
 R_tR_s=\rho(\gamma)R_r z_{t,s},
\]
where all $z_{t,s}$ lie in the image $Z_c$ of the finite central
kernel of the covering map. For $C_t=B_tR_t^{-1}$ one has $C_t\rho(g)=\rho(g)C_t$ for every $g\in G_0$. Hence
\[
 C_tC_s\in K C_r Z_c,\qquad C_tKC_t^{-1}=K.
\]
It follows, also using inverses, that
\[
 \mathcal C=\bigcup_{t\in\mathcal T}K Z_c C_t
\]
 is a compact group of analytic germs fixing zero. It is a finite
union of compact $K$-cosets. The finitely many representatives and
compactness of $K$ give common domains for the averaging below.
Thus $c\rho(g)=\rho(g)c$ for $(c,g)\in\mathcal C\times G_0$.

With $m_{\mathcal C}$ the Haar probability on $\mathcal C$, the Bochner average~\cite{Bochner} is the analytic germ
\[
 B(w)=\int_{\mathcal C}(D_0c)^{-1}c(w)\,dm_{\mathcal C}(c),
 \qquad D_0B=I,\qquad Bc=(D_0c)B.
\]
The integral converges on a common complex ball. For $a\in\mathcal C$ substitute $c\mapsto ca$ and use $D_0(ca)=D_0c\,D_0a$ to obtain $B(aw)=D_0a\,B(w)$. The inverse function theorem applies since $D_0B=I$. Every integrand commutes with $\rho(G_0)$, and $Bk=(D_0k)B=kB$ for $k\in K$.
Replace $\psi_p$ by $\psi_pB^{-1}$. Then
\[
 t\psi_p=\psi_{T_tp}A_t,\qquad A_t=(D_0C_t)R_t\in\GL(V),
 \qquad A_tKA_t^{-1}=K.
\]
Here $BkB^{-1}=k\in\mathrm{GL}(V)$ for $k\in K$, and the compact image $\{D_0c:c\in\mathcal C\}$ satisfies $(D_0c)K(D_0c)^{-1}=K$.

\par\Needspace{3\baselineskip}\smallskip\noindent\textbf{Analyticity of the lift.} Identify a frame with its first derivative $b=D_0\psi_p:V\to T_xM$, where $x=c(p)$. The candidate first frame over $tx$ is $D_xt\,bA_t^{-1}$. In each analytic section from Lemma~\ref{sections}, membership in the frame bundle is a closed condition since its fiber in these coordinates is $K$. Formula~\eqref{H-eq:finite-lift} puts this frame in the bundle for $\eta$-almost every $p\in P_1$ and hence for every $p\in P_1$ by continuity and full support. The germ identity extends there coefficient by coefficient. The inverse bundle map is given by the inverse first-frame formula over $t^{-1}x$. For all $p\in P_1$ and $h\in K$, the derivative formula gives
\[
 c(T_t p)=t c(p),\qquad T_t(ph)=T_t(p)(A_t h A_t^{-1}).
\]
The equality $T_tT_s(p)=T_{ts}(p)k_{t,s}$ extends from a set of full $\eta$-measure by continuity. The analytic lift agrees $\eta$-almost everywhere with the measured lift in \eqref{H-eq:finite-lift}, hence preserves $\eta$.
\end{proof}
\section{The coframe from second derivatives}\label{sec:coframe-from-jets}

Write $P_1\to O_1'$ for the bundle of Lemmas~\ref{sections} and~\ref{finite}, $c(p)=\psi_p(0)$ for its projection, and $\mathfrak k=\Lie K$. Its frame probability is smooth and positive, since $m_1|_{O_1}=k\,m|_{O_1}$, and the measured $\Lambda$-action is ergodic by Proposition~\ref{action-coordinates}. Use the source coordinate of Lemma~\ref{finite}, so that $t\psi_p=\psi_{T_tp}A_t$ and $A_tKA_t^{-1}=K$.

\begin{lemma}[The coframe and its covariance]\label{coframe}
There is an analytic family of linear isomorphisms $\vartheta_p:T_pP_1\to V\oplus\mathfrak k$ with $\vartheta_{\gamma p}D_p\gamma=(\rho(\gamma)\oplus\id)\vartheta_p$ for $\gamma\in\Lambda$. For $t\in\Gamma_1$ the lift from Lemma~\ref{finite} satisfies $T_t^*\vartheta=(A_t\oplus\Ad(A_t)|_{\mathfrak k})\vartheta$. For $z\in\mathfrak k$ its value on $z^\#(p)=\partial_s|_0\,p\exp(sz)$ is $(0,z)$.
\end{lemma}
\begin{proof}
Fix analytic coordinates on an open subset of $O'_1$ and write $x=c(p)$ in these coordinates. The matrix $b(p)=D_0\psi_p$ depends analytically on $p$ by Lemma~\ref{sections}. For $\xi\in T_pP_1$, $db(\xi)$ is its derivative along any curve of frames with initial velocity $\xi$. Write $b=b(p)$. For $v,w$ in the ambient coordinate space set
\[
 C_x(v,w)=-D_0^2\psi_p(b^{-1}v,b^{-1}w).
\]
Replacing $p$ by $ph$ replaces $b$ by $bh$ and $D_0^2\psi_p$ by $D_0^2\psi_p\,(h\otimes h)$, leaving $C_x$ unchanged. By complete reducibility choose a $\rho(G_0)$-equivariant projection $\pi_0:\operatorname{End}(V)\to\mathfrak k$ for conjugation. In Lemma~\ref{finite}, $D_0\mathcal C=\{D_0a:a\in\mathcal C\}$ is compact, contains $K$, normalizes $\mathfrak k$ and centralizes $\rho(G_0)$. With $dh$ its Haar probability, set
\[
 \operatorname{pr}(A)=\int_{D_0\mathcal C} h^{-1}\pi_0(hAh^{-1})h\,dh
 \qquad (A\in\operatorname{End}(V)).
\]
It is the identity on $\mathfrak k$. Right invariance of $dh$ gives $\operatorname{pr}(hAh^{-1})=h\operatorname{pr}(A)h^{-1}$ for $h\in D_0\mathcal C$. Commutation with $\rho(G_0)$ gives the same identity for $h\in\rho(G_0)$ and hence for each $A_t\in (D_0\mathcal C)\rho(G_0)$. For $\xi\in T_pP_1$ put $dx(\xi)=D_pc(\xi)$ in the ambient chart and set
\[
 \vartheta_p(\xi)=\left(b^{-1}dx(\xi), \operatorname{pr}\!\left[b^{-1}db(\xi)+b^{-1}C_x(dx(\xi),b\,\cdot)\right]\right).
\]
Under a change of target coordinates $x'=a(x)$ we have $b'=D_xa\,b$ and
\[
 \begin{aligned}
 C'_{x'}(D_xa\,v,D_xa\,w)&=D_xa\,C_x(v,w)-D_x^2a(v,w),\\
 db'(\xi)&=D_xa\,db(\xi)+D_x^2a(dx(\xi),b\,\cdot).
 \end{aligned}
\]
Substitute $v=dx(\xi)$ and $w=b\,\cdot$ in the first line, add the second line, and multiply by $(b')^{-1}=b^{-1}(D_xa)^{-1}$. This gives exactly the expression inside $\operatorname{pr}$ in the original coordinates. On a vertical tangent $z^\#$ one has $dx=0$ and $db=bz$, giving $(0,z)$. The first component identifies $T_pP_1/\ker D_pc$ with $V$, so $\vartheta_p$ is an isomorphism. For a constant source change $\psi_p\mapsto\psi_p\circ A^{-1}$ we have $b'=bA^{-1}$ and $C'_x=C_x$. Hence the first coframe entry becomes $A b^{-1}dx(\xi)$, and the expression inside $\operatorname{pr}$ is conjugated by $A$. For $A=\rho(\gamma),A_t,h^{-1}$ the projection commutes with this conjugation. Substitution in \eqref{frames} and Lemma~\ref{finite} gives the two asserted covariance identities and $R_h^*\vartheta=(h^{-1}\oplus\operatorname{Ad}(h^{-1}))\vartheta$. Here $\rho(\gamma)$ centralizes $\mathfrak k$, so it acts trivially on the second entry.
\end{proof}

We identify the one-parameter changes of Part~\ref{agreement} in this coframe. Let $v_t$ be one of the families in Lemma~\ref{lie} or its $G_0$-conjugates in Section~\ref{sec:analytic-sections}. Recall $H_{[g,p]}=\iota_g\psi_p\rho(g)^{-1}$ from \eqref{susp}. By Lemma~\ref{joining}, $H_{T_{v_t}z}=H_z\circ v_t$ for $\widehat m$-almost every $z\in\widehat X$ at each fixed small $t$. Disintegrate over $G_0/\Lambda$ and intersect over the rational $t$. Fix a coset $g\Lambda$ in the resulting set of full Haar measure and a representative $g\in G_0$. In the coordinate $\psi_p B^{-1}$ of Lemma~\ref{finite}, formula~\eqref{susp} changes $v_t$ to $u_t=B\rho(g)^{-1}v_t\rho(g)B^{-1}$. For $w$ near $0$ in $V$ put $Y(w)=\partial_t|_0 u_t(w)$. With $\psi_p$ now denoting the chart in this coordinate, the moved frame $p_t$ has center $\psi_p(u_t(0))$ and first derivative $D_0(\psi_p\circ u_t)$. These formulas are analytic on a small frame neighborhood for small $t$. For each small rational $t$, the candidate belongs to the bundle and $\psi_{p_t}=\psi_p\circ u_t$ for $\eta$-almost every $p\in P_1$ in this neighborhood. The bundle is locally closed over $O_1'$. Full support and continuity in $p$, then continuity in $t$, give both statements throughout the neighborhood, applying continuity to each Taylor coefficient for the germ identity.

With $b=D_0\psi_p$ and dots denoting derivatives at $t=0$, the chain rule gives
\[
 \dot x=bY(0),\qquad \dot b=D_0^2\psi_p(Y(0),\cdot)+bD_0Y,\qquad
 \vartheta_p(\dot p_0)=\bigl(Y(0),\operatorname{pr}(D_0Y)\bigr).
\]
The second derivative cancels in the formula of Lemma~\ref{coframe}. Since $u_{t+s}=u_t\circ u_s$, the same coframe value gives the velocity at $p_t$. Thus $p_t$ is the local flow, through $p$, of the field $\vartheta^{-1}(Y(0),\operatorname{pr}(D_0Y))$. Its value in the coframe is independent of $p$.

\section{Integrating the coframe}
Let $\mathfrak l=V\oplus\mathfrak k$ be the value space of the coframe of Lemma~\ref{coframe}. We first determine its Lie bracket, then integrate its fields. For $u\in\mathfrak l$, $p\in P_1$ and $a\in G_0$ put
\[
 X_u(p)=\vartheta_p^{-1}u,\qquad
 \varrho(a)(v,z)=(\rho(a)v,z)\quad(v\in V,\ z\in\mathfrak k).
\]
Thus $D_p\gamma X_u(p)=X_{\varrho(\gamma)u}(\gamma p)$ for $\gamma\in\Lambda$. Recall from Lemma~\ref{sections} that $\eta$ is smooth and positive on $P_1$ and that $\Lambda$ acts ergodically.
\begin{lemma}[Coframe volume]\label{H-thm:density}
For a Lebesgue density $du$ on $\mathfrak l$ put $v_\vartheta=\vartheta^*du$. Then $\eta=v_\vartheta/v_\vartheta(P_1)$, and $P_1$ has finitely many connected components.
\end{lemma}
\begin{proof}
The character $\det\varrho$ of the connected semisimple group $G_0$ is trivial. Covariance of $\vartheta$ gives $\gamma^*v_\vartheta=v_\vartheta$ for $\gamma\in\Lambda$. The positive smooth function $d\eta/dv_\vartheta$ is therefore invariant and constant by ergodicity. Integrating over $P_1$ gives its value. If $D$ is a connected component, each distinct translate has the same mass $\eta(D)>0$. Thus there are at most $1/\eta(D)$ such translates. Their union is invariant of positive measure, and its complement is an open null set. Full support gives $P_1=\bigcup_{\gamma\in\Lambda}\gamma D$.
\end{proof}

\subsection*{Invariant fields}
Let $T:P\to P$ be a smooth diffeomorphism of a second countable smooth manifold, and let $\eta\sim\vol_P$ be a $T$-invariant ergodic probability. For a smooth field $X:P\to TP$, let $I_p\ni0$ denote the maximal time interval of the trajectory from $p$. For $t\in\mathbb R$, put $P_t=\{p\in P:t\in I_p\}$. The smooth diffeomorphism $\phi^t:P_t\to P_{-t}$ satisfies $\partial_t\phi^t(p)=X(\phi^t(p))$, $\phi^0(p)=p$, and $(\phi^t)^{-1}=\phi^{-t}$.

\begin{lemma}[Completeness of an invariant field]\label{H-lem:flow}
If $D_pT\,X(p)=X(Tp)$ for every $p\in P$, then $\eta\{p\in P:I_p=\mathbb R\}=1$ and $(\phi^t)_*\eta=\eta$ for every $t\in\mathbb R$.
\end{lemma}
\begin{proof}
Write $(-\tau_-(p),\tau_+(p))$ for the maximal time interval
through $p$. The measurable functions $\tau_\pm:P\to(0,\infty]$
satisfy
\[
 T\phi^t p=\phi^t Tp,\qquad \tau_\pm(Tp)=\tau_\pm(p).
\]
Ergodicity gives constants $a_\pm\in(0,\infty]$ with $\tau_\pm=a_\pm$ $\eta$-almost everywhere. For $|t|<\min(a_+,a_-)$ both flow domains have full $\eta$-measure. Since $\phi^t:P_t\to P_{-t}$ is a diffeomorphism and $\eta\sim\vol_P$, the probability $\eta_t=(\phi^t)_*(\eta|_{P_t})$ is equivalent to $\eta$. Commutation gives $T_*\eta_t=\eta_t$, so $d\eta_t/d\eta$ is $T$-invariant. Ergodicity and $\eta_t(P)=1$ give $\eta_t=\eta$.
If $a_+<\infty$, choose such a $t\in(0,a_+)$. The lifetime
identity and invariance give
\[
 a_+=\tau_+(\phi^t p)=\tau_+(p)-t=a_+-t
 \quad(\eta\text{-a.e. }p\in P),
\]
a contradiction. Apply the same argument to $-X$.
Now $a_+=a_-=\infty$, and iteration of small times proves flow
invariance for every $t$.
\end{proof}

\subsection*{The coframe action on the measure algebra}

\begin{definition}[Actions on the measure algebra]\label{H-def:probability-algebra-action}
For a probability space $(P,\eta_P)$, its measure algebra $\mathcal A$ consists of the classes $[A]$ of measurable $A\subset P$ under $A\sim B\Longleftrightarrow\eta_P(A\mathbin\triangle B)=0$, with $d([A],[B])=\eta_P(A\mathbin\triangle B)$. A probability-preserving action by Boolean automorphisms $R_g$ satisfies $\eta_P(R_g[A])=\eta_P(A)$ for $g\in\mathcal G$, $[A]\in\mathcal A$. For a left action $R_{gh}=R_gR_h$, and for a right action $R_{gh}=R_hR_g$. Writing $gA$ for a representative of $R_g[A]$, the action is strongly continuous when
$\eta_P(gA\mathbin{\triangle}A)\to0$ as $g\to e$ for every
measurable $A$. \end{definition}
For indicator functions, $\|{\bf1}_{gA}-{\bf1}_A\|_2^2=\eta_P(gA\mathbin\triangle A)$. Finite linear combinations of indicators are dense in $L^2(P,\eta_P)$, and the Koopman operators are isometries. Hence Definition~\ref{H-def:probability-algebra-action} gives $\lim_{g\to e}\|U_g f-f\|_{L^2(\eta_P)}=0$ for all $f\in L^2(\eta_P)$, and taking $f=\mathbf1_A$ gives the converse.

Return to $X_u(p)=\vartheta_p^{-1}u$ on $P_1$, where $\mathfrak l=V\oplus\mathfrak k$. For fixed $u\in\mathfrak l$ and $t\in\mathbb R$ write
\[
 P_{u,t}=\{p\in P_1:t\text{ belongs to the maximal existence interval of }X_u\text{ at }p\},
 \qquad \phi_u^t:P_{u,t}\longrightarrow P_{u,-t}.
\]
Write $\mathcal A$ for the measure algebra of $(P_1,\eta)$ as in Definition~\ref{H-def:probability-algebra-action}.

\Needspace{10\baselineskip}
\begin{theorem}[Integration of the coframe]\label{H-thm:coframe-integration}
There is a Lie bracket on $\mathfrak l$ such that $[X_u,X_v]=X_{[u,v]}$ for $u,v\in \mathfrak l$. For the simply connected group $\widetilde L$ of this Lie algebra there is a strongly continuous measure-preserving right action $g\mapsto R_g$ on $\mathcal A$, $R_{gh}=R_hR_g$. For each fixed $u\in \mathfrak l$ and $t\in\mathbb R$,
\[
  R_{\exp(tu)}[A]=[\phi_u^t(A)]
 \qquad(A\subset P_{u,t}\text{ Borel}).
\]
For $\gamma\in\Lambda$, the automorphism $\sigma_\gamma$ of $\widetilde L$ with $D_e\sigma_\gamma=\varrho(\gamma)$ satisfies $\gamma R_g\gamma^{-1}=R_{\sigma_\gamma(g)}$ on $\mathcal A$.
\end{theorem}

\begin{proof}
\textbf{The Lie bracket.} For $u,v\in\mathfrak l$ put $C(p)(u,v)=\vartheta_p([X_u,X_v]_p)$. Covariance gives
\[
 C(\gamma p)(\varrho(\gamma)u,\varrho(\gamma)v)
                  =\varrho(\gamma)C(p)(u,v).
\]
For $g\cdot A=\varrho(g)A(\varrho(g)^{-1}\cdot,\varrho(g)^{-1}\cdot)$ on $\Hom(\bigwedge^2\mathfrak l,\mathfrak l)$, covariance and $\gamma_*\eta=\eta$ give $(A\mapsto\gamma\cdot A)_*(C_*\eta)=C_*\eta$ for all $\gamma\in\Lambda$. Lemma~\ref{H-lem:fixed-law} and ergodicity give $C_*\eta=\delta_{C_0}$ with $C_0\in\Hom_{G_0}(\bigwedge^2 \mathfrak l,\mathfrak l)$. Since $\operatorname{supp}\eta=P_1$ and $p\mapsto C(p)$ is continuous, $C(p)=C_0$ for all $p\in P_1$. Because $[X_u,X_v]=X_{C_0(u,v)}$, applying $\vartheta$ to the Jacobi identity of these fields yields $\sum_{\rm cyc} C_0(u,C_0(v,w))=0$ for $u,v,w\in \mathfrak l$. Equip $\mathfrak l$ with the bracket $[u,v]=C_0(u,v)$. Let $\widetilde L$ be its simply connected Lie group, and let $u^{\widetilde L}$ be the left-invariant field with value $u$ at the identity. The fields $(u^{\widetilde L},X_u)$ satisfy
\[
 [(u^{\widetilde L},X_u),(v^{\widetilde L},X_v)]
       =([u,v]^{\widetilde L},X_{[u,v]}).
\]
Projection on the first factor identifies their span with $T\widetilde L$. Frobenius gives an integral leaf through $(e,p)\in\widetilde L\times P_1$. The derivative of its first projection is an isomorphism, so the inverse function theorem writes the leaf locally as $g\mapsto(g,p\cdot g)$ for $g$ near $e$ in $\widetilde L$. Left translation on the first factor preserves $(u^{\widetilde L},X_u)$. Uniqueness of the integral leaf therefore gives $(p\cdot g)\cdot h=p\cdot(gh)$ whenever these local products are defined.

\par\Needspace{3\baselineskip}\smallskip\noindent\textbf{Completeness in weight directions.} Form the suspension $D=(G_0\times P_1)/\Lambda$ with $(g,p)\sim(g\gamma^{-1},\gamma p)$, projection $[g,p]\mapsto g\Lambda$, and probability $\eta_D$ induced by Haar probability on $G_0/\Lambda$ and $\eta$ on the fibers. Its $G_0$-action is ergodic by the suspension argument after \eqref{susp}, since an invariant function is constant in $g$ and then $\Lambda$-invariant in $p$. Its Koopman representation is strongly continuous by density of compactly supported continuous functions. For $z=[g,p]\in D$ and $u\in\mathfrak l$, let $Y_u(z)$ be the tangent represented at $(g,p)$ by $(0,X_{\varrho(g)^{-1}u}(p))$. Covariance $D\gamma X_v=X_{\varrho(\gamma)v}\circ\gamma$ gives independence of the representative and $a_*Y_u=Y_{\varrho(a)u}$. Taking brackets at fixed $g$ gives
\[
 [Y_u,Y_v]_{[g,p]}=[\,0,[X_{\varrho(g)^{-1}u},X_{\varrho(g)^{-1}v}]_p\,]
                 =Y_{[u,v]}([g,p]).
\]
Thus on $\widetilde L\times D$ the distribution of $(u^{\widetilde L},Y_u)$ is involutive. Its projection onto $T\widetilde L$ is an isomorphism on each tangent plane. Frobenius gives the local group words on the suspension used below, on one parameter neighborhood after restricting the starting point to a compact subchart.

Decompose $\mathfrak l=\bigoplus_\chi\mathfrak l_\chi$ with $\mathfrak l_\chi=\{u\in\mathfrak l:\varrho(\exp H)u=e^{\chi(H)}u\ (H\in\mathfrak a)\}$. Choose a basis $u_1,\ldots,u_n$ of weight vectors, $n=\dim\mathfrak l$. Thus $\varrho(\exp H)u_i=e^{\chi_i(H)}u_i$ for $H\in\mathfrak a$. By Lemma~\ref{H-lem:hyperplane}, choose $H_i\in\ker\chi_i$ with nonzero projection on every simple factor, so $\exp H_i$ is ergodic on $(D,\eta_D)$. Since $\varrho(\exp H_i)u_i=u_i$, the covariance above gives $(\exp H_i)_*Y_{u_i}=Y_{u_i}$. Write $\phi_i^t$ for the flow of $Y_{u_i}$ on the suspension $D$. Applying Lemma~\ref{H-lem:flow} with $T(z)=(\exp H_i)z$ and $X=Y_{u_i}$ gives $\eta_D\{z\in D:\phi_i^t(z)\text{ exists for all }t\in\mathbb R\}=1$ and $(\phi_i^t)_*\eta_D=\eta_D$ for every $t$. Its Koopman group on $L^2(D,\eta_D)$, $U_i(t)f=f\circ\phi_i^t$, satisfies $\lim_{t\to0}\|U_i(t)f-f\|_{L^2(\eta_D)}=0$, first for bounded continuous $f$ by bounded convergence, then for all $f\in L^2(\eta_D)$ by density and unitarity.

\par\Needspace{3\baselineskip}\smallskip\noindent\textbf{Relations between pairs of flows.} Fix $1\le i,j\le n$. If some $a=\exp H$, $H\in\mathfrak a$, satisfies $\chi_i(H)<0$ and $\chi_j(H)<0$, use $\bigoplus_{\chi(H)<0}\mathfrak l_\chi$. This is a subalgebra since $[\mathfrak l_\chi,\mathfrak l_{\chi'}]\subset \mathfrak l_{\chi+\chi'}$. Otherwise $\dim\operatorname{span}\{\chi_i,\chi_j\}\le1$, since two independent forms take the pair of values $(-1,-1)$ somewhere. Apply Lemma~\ref{H-lem:hyperplane} to a linear form spanning that line (zero if both weights vanish). The resulting ergodic element $b$ fixes both fields. Use the subalgebra fixed by $b$.
In either subalgebra, use the subset of the chosen weight basis that belongs to it, in its original order. The local group coordinates are
$w(c)=\prod_k\exp(c_k u_k)$ near $0$ in the real coordinate space indexed by these $k$, in basis order. For all sufficiently
small real $s,t$ define $c(s,t)$ by
\[
  w(c(s,t))=\exp(-tu_i)\exp(su_j)\exp(tu_i).
\]
The inverse function theorem makes $c$ analytic, $c(0,t)=0$,
and $\sum_k\partial_s c_k(0,t)u_k=\Ad(\exp(-tu_i))u_j$.
\begin{claim}
For every sufficiently small fixed pair $(s,t)\in\mathbb R^2$, the operators $U_i(t)f=f\circ\phi_i^t$ on $L^2(D,\eta_D)$ satisfy
\begin{equation}
 U_j(s)U_i(t)=U_i(t)\prod_k U_k(c_k(s,t)).
 \label{H-eq:coframe-pair}
\end{equation}
\end{claim}
\begin{claimproof}
Under common contraction conjugate both words by $a^n$. Equivariance
of the ordered word coordinates gives
\[
 c_k(e^{n\chi_j(\log a)}s,e^{n\chi_i(\log a)}t)
       =e^{n\chi_k(\log a)}c_k(s,t).
\]
Here $\varrho(a)u_k=e^{\chi_k(\log a)}u_k$. For the negative weights in this word, $e^{n\chi_k(\log a)}\to0$ as $n\to\infty$. Cover the suspension by countably
many relatively compact flow patches $U$. Recurrence gives $\#\{n\ge1:a^nx\in U\}=\infty$ for $\eta_D|_U$-almost every $x\in U$ simultaneously over this countable cover.
For $\eta_D$-almost every starting $x\in D$, choose such a patch $U\ni x$ with returns $a^nx\in U$ for arbitrarily large $n$. For every sufficiently large return $n$ all the parameters $e^{n\chi_k(\log a)}c_k$ lie in its common word domain. The covariance $a^n\phi_k^s a^{-n}=\phi_k^{e^{n\chi_k(\log a)}s}$ transports the local identity at $a^nx$ back to the required identity at $x$. For these fixed $(s,t)$ the union of the exceptional sets over $n\ge1$ has $\eta_D$-measure zero. In the fixed subalgebra $b_*Y_{u_k}=Y_{u_k}$ for every field in the two words, hence $b\phi_k^t=\phi_k^t b$. For each small fixed $(s,t)$, let $A_{s,t}\subset D$ consist of the starting points $z$ where both flow words in \eqref{H-eq:coframe-pair} are defined and have the same endpoint. The local action gives an open patch $U\subset A_{s,t}$ with $\eta_D(U)>0$. It follows that $\eta_D(A_{s,t}\mathbin\triangle b^{-1}A_{s,t})=0$ and $\eta_D(A_{s,t})>0$. Ergodicity gives $\eta_D(A_{s,t})=1$. This proves \eqref{H-eq:coframe-pair}.
\end{claimproof}
\textbf{Integration of the pairwise relations.}
For $f\in L^2(D,\eta_D)$ let $A_i f=\lim_{t\to0}(U_i(t)f-f)/t$ whenever the limit exists in $L^2(D,\eta_D)$. Put $\mathcal D=\bigcap_i\operatorname{Dom}A_i$. We show $C_c^\infty(D)\subset\mathcal D$, hence $\overline{\mathcal D}^{\,L^2(\eta_D)}=L^2(\eta_D)$. For $f\in C_c^\infty(D)$ integrate $Y_{u_i}f$ along the
trajectories defined for all $t$ on an $\eta_D$-full set. The identity $U_i(t)f-f=\int_0^t U_i(s)Y_{u_i}f\,ds$ and
strong continuity show $A_i f=Y_{u_i}f$ in $L^2$. For $v=\sum v_k u_k$
write $A_v:\mathcal D\to L^2(D,\eta_D)$, $A_v f=\sum_k v_k A_k f$.
Differentiate \eqref{H-eq:coframe-pair} at $s=0$ on $f\in\mathcal D$. Indeed let $D_k(s)=U_k(c_k(s,t))$. Since $D_k(0)=I$,
\[
 \frac{(\prod_kD_k(s)-I)f}{s}
 =\sum_k\left(\prod_{l<k}D_l(s)\right)
                  \frac{(D_k(s)-I)f}{s}
 \longrightarrow\sum_k\partial_s c_k(0,t)A_k f
 \quad\text{in }L^2 .
\]
The finite products are unitary and strongly continuous. On the left of \eqref{H-eq:coframe-pair}, existence of this derivative is exactly $U_i(t)f\in\operatorname{Dom}A_j$. Thus $U_i(t)f\in\mathcal D$ and
\begin{equation}
 A_v U_i(t)f
    =U_i(t)A_{\Ad(\exp(-tu_i))v}f .
 \label{H-eq:coframe-adj}
\end{equation}
Linearity gives \eqref{H-eq:coframe-adj} for $v\in \mathfrak l$. For $t\in\mathbb R$ apply it $m$ times with $t/m$ sufficiently small. Applying domain preservation also to $-t$ gives $U_i(t)\mathcal D=\mathcal D$. Formula \eqref{H-eq:coframe-adj} and strong continuity also give continuity
in the norm $\|f\|_2+\sum_i\|A_i f\|_2$ on this domain.

On an identity neighborhood $O\subset\widetilde L$ with unique ordered coordinates $g=\prod_i\exp(t_i u_i)$, put $F:O\to\mathcal U(L^2(D,\eta_D))$, $F(g)=\prod_iU_i(t_i)$. Each factor is unitary and strongly continuous in $t_i$.
\begin{claim}
$F(gh)=F(g)F(h)$ for $g,h$ in a sufficiently small identity neighborhood.
\end{claim}
\begin{claimproof}
For a differentiable path $g(s)$ in the coordinate chart and $f\in\mathcal D$, \eqref{H-eq:coframe-adj} gives
\begin{equation}
 \frac{d}{ds}F(g(s))f=F(g(s)) A_{g(s)^{-1}g'(s)}f .
 \label{H-eq:coframe-product}
\end{equation}
Indeed, \eqref{H-eq:coframe-adj} moves the $i$th generator past the word on its right, conjugating $t'_iu_i$ by the inverse of that word. Explicitly, with $g_i=\exp(t_i u_i)$ the sum is $\sum_i\operatorname{Ad}((g_{i+1}\cdots g_n)^{-1})(t_i'u_i)=g^{-1}g'$ by the product rule in the group. Differentiation by telescoping is valid on $\mathcal D$ since every partial word preserves it. The inverse difference identity then gives $(F(g)^{-1})'f=-A_{g^{-1}g'}F(g)^{-1}f$. For fixed $g$ and a path $h(s)$ from $e$ to $h$ in a smaller coordinate ball, $(gh)^{-1}(gh)'=h^{-1}h'$. Hence on $\mathcal D$
\[
 \frac d{ds}\bigl(F(gh(s))F(h(s))^{-1}f\bigr)=0 .
\]
For the product rule take the first difference on the fixed vector $F(h(s))^{-1}f$ and the second using strong continuity and the unitary bound. Integrating from $h(0)=e$ to $h$ and using density of $\mathcal D$ gives the claimed identity.
\end{claimproof}

\par\Needspace{3\baselineskip}\smallskip\noindent\textbf{Descent to the frame bundle.}
The flows on $D=(G_0\times P_1)/\Lambda$ preserve the base coset. Let $m_{G_0/\Lambda}$ be Haar probability and $\eta_{g\Lambda}=(p\mapsto[g,p])_*\eta$ the conditional on its fiber. For rational $t$, invariance of $\eta_D=\int\eta_{g\Lambda}\,dm_{G_0/\Lambda}(g\Lambda)$ and uniqueness of disintegration give $(\phi_i^t)_*\eta_{g\Lambda}=\eta_{g\Lambda}$ for $m_{G_0/\Lambda}$-almost every $g\Lambda\in G_0/\Lambda$. Disintegrate also the full-measure set on which the finitely many basis fields have infinite existence intervals.

Take a countable dense set of pairs $(a,b)$ in the local multiplication domain. For each pair, the identity $F(ab)=F(a)F(b)$ on the suspension is equality of measurable transformations modulo $\eta_D$ (test on a countable family of Borel sets separating points). Intersect the full-measure sets for these pairs and the rational times, and fix one coset $g_0\Lambda$ in the intersection. Continuity of the complete trajectories and bounded convergence imply invariance of $\eta_{g_0\Lambda}$ for every real time and strong continuity of each Koopman group. Therefore on $L^2(\eta_{g_0\Lambda})$,
\[
 F(a_nb_n)=F(a_n)F(b_n),\quad (a_n,b_n)\longrightarrow(a,b)
 \quad\Longrightarrow\quad F(ab)=F(a)F(b).
\]
The operator equality now holds for every pair in the local domain.

For a path $a:[0,1]\to\widetilde L$ from $e$ to $g$, choose a subdivision $0=t_0<\cdots<t_n=1$ so small that all increments and their partial increments lie in the local domain. Assign to it
\[
 \prod_{j=1}^n F\big(a(t_{j-1})^{-1}a(t_j)\big).
\]
The local multiplication law gives the same product after refinement. Subdivide a homotopy with endpoints fixed into rectangles small enough for the same law to hold along either route around each rectangle. The products on its two boundary paths agree. Since $\widetilde L$ is simply connected, the product depends only on $g$. Concatenation gives the group law. Each operator is a product of Koopman operators, and hence gives an automorphism of the measure algebra.

In the fiber over $g_0\Lambda$ the field with parameter $u_i$ is $X_{v_i}$ on $P_1$ with $v_i=\varrho(g_0)^{-1}u_i$. Thus Lemma~\ref{H-lem:flow} on the suspended weight fields gives
\begin{equation}\label{H-eq:complete-basis}
 \eta\{p\in P_1:\phi_{v_i}^t(p)\text{ exists for all }t\in\mathbb R\}=1
 \qquad (1\le i\le\dim \mathfrak l).
\end{equation}
In this fiber the group parameter $\exp(tu)$ corresponds on $P_1$ to $X_{\varrho(g_0)^{-1}u}$. Substituting $u=\varrho(g_0)v$ gives $X_v$. Thus precompose the fiber action with the unique automorphism of the simply connected group whose differential is $\varrho(g_0)$.
On a relatively compact open set in $P_1$ where the local action is defined, the ordered product of basis exponentials for $\exp(tv)$
equals the local flow of $X_v$. For each fixed word parameter,
the complete factors agree with their analytic local flows for $\eta$-almost every starting frame $p$ in the domain. For each local word with parameter $b\in\mathbb R^n$, Fubini gives agreement for $(\eta\otimes\mathcal L^n)$-almost every $(p,b)$ in its domain,
as required in Lemma~\ref{H-lem:coframe-measured-realization}. Each factor preserves $\eta$. A countable chart cover and time subdivision therefore extend the identity, for each fixed time $t$, to $\eta$-almost every $p\in P_1$ for which the analytic trajectory reaches $t$. For a fixed $\gamma\in\Lambda$, $D\gamma X_{v_i}=X_{\varrho(\gamma)v_i}\circ\gamma$. Since $\gamma_*\eta=\eta$, uniqueness of trajectories and \eqref{H-eq:complete-basis} give $\gamma R_{\exp(tv_i)}\gamma^{-1}=R_{\exp(t\varrho(\gamma)v_i)}$ for each $t$ on the measure algebra ($[A]=[B]\iff\eta(A\mathbin\triangle B)=0$, Definition~\ref{H-def:probability-algebra-action}). Products of these basis exponentials contain an identity neighborhood and generate the connected group $\widetilde L$. This gives the asserted conjugacy for all $g$.
\end{proof}

\section{From frame orbits to an affine action}
Return to the frame bundle $P_1\to O_1'$ of Lemmas~\ref{finite} and~\ref{coframe}, with group $K$, probability $\eta$ and coframe $\vartheta$. Theorem~\ref{H-thm:coframe-integration} gives the measured action of the simply connected group $\widetilde L$ with coframe Lie algebra $\mathfrak l=V\oplus\mathfrak k$. We first include the disconnected part of $K$ in the frame action.

\begin{lemma}[The action of the full frame group]\label{adjoin}
There is a Lie group $L$ acting faithfully on the measure algebra of $(P_1,\eta)$, containing $K$ with $R_h(p)=ph$. It has the coframe Lie algebra and finitely many connected components. Each $T_t$ in Lemma~\ref{finite} normalizes this action by an automorphism preserving $K$.
\end{lemma}
\begin{proof}
By Lemmas~\ref{finite} and~\ref{coframe}, $T^*\vartheta=B\vartheta$ with
\[
 B=A_t\oplus\Ad(A_t)|_{\mathfrak k}\quad(T=T_t),\qquad
 B=h^{-1}\oplus\Ad(h^{-1})|_{\mathfrak k}\quad(T=R_h,\ h\in K).
\] Then $T_*X_u=X_{Bu}$ and $B[u,v]=[Bu,Bv]$. The simply connected group integrates $B$ to an automorphism. By \eqref{H-eq:complete-basis} there is a basis $v_i$ with complete trajectories for $\eta$-almost every starting point. Applying $T$ to those trajectories and using uniqueness gives $TR_{\exp(sv_i)}T^{-1}=R_{\exp(sBv_i)}$ modulo $\eta$ for each $s$. Products of these exponentials generate $\widetilde L$.

Thus the semidirect product $\widetilde L\rtimes K$ acts on the measure algebra, where $h$ acts on $\operatorname{Lie}(\widetilde L)$ by $(v,z)\mapsto(hv,\Ad(h)z)$. Let $N$ be the kernel of this action on the measure algebra and put $L=(\widetilde L\rtimes K)/N$. Strong continuity makes $N$ closed. For $h\in K$, if $R_h[A]=[A]$ for every Borel $A\subset P_1$, a countable family of Borel sets separating points gives $ph=p$ for $\eta$-almost every $p\in P_1$. Freeness gives $h=e$, so $K\to L$ is injective. For $z\in\mathfrak k$, put $\iota(z)=(0,z)$. The coframe value of the vertical derivative of $p\exp(tz)$ is $\iota(z)$. Since $\Ad(\exp_K(tz))z=z$, the pair
\[
 t\longmapsto\bigl(\exp_{\widetilde L}(-t\iota(z)),\,\exp_K(tz)\bigr)
\]
is a one-parameter subgroup of $\widetilde L\rtimes K$. Its two factors act on $P_1$ by the vertical flows $p\mapsto p\exp_K(-tz)$ and $p\mapsto p\exp_K(tz)$. Their product is the identity, so $(-\iota(z),z)\in\operatorname{Lie}(N)$. Conversely a kernel vector $(u,z)$ acts trivially for every rational time on a common set of full $\eta$-measure. In a local flow chart differentiation at zero gives $\vartheta_p^{-1}(u+\iota(z))=0$ there, so $u=-\iota(z)$. Therefore $\operatorname{Lie}(N)=\{(-\iota(z),z):z\in\mathfrak k\}\subset\mathfrak l\oplus\mathfrak k$. This gives the asserted Lie algebra and $|\pi_0(L)|\le|\pi_0(K)|$. Conjugation by $T_t$ on the two factors is the integrated automorphism above and $h\mapsto A_thA_t^{-1}$, respectively. It preserves the kernel, so descends to $L$.
\end{proof}

\Needspace{7\baselineskip}
\begin{lemma}[Injectivity from a set of full measure]\label{H-lem:conull-local-diffeomorphism}
Let $j:U\to M$ be a local diffeomorphism between second countable
smooth manifolds of the same dimension. Let $\xi\sim\vol_U$ and let $C\subset U$ be Borel with
\[
                 \xi(U\setminus C)=0,\qquad u,v\in C,\ j(u)=j(v)\ \Longrightarrow\ u=v.
\]
Then $j$ is a diffeomorphism onto the open set $j(U)$.
\end{lemma}
\begin{proof}
If $u\ne v$ and $j(u)=j(v)$, choose disjoint neighborhoods $V_u,V_v$
on which $j$ is a diffeomorphism onto the same nonempty open set
$W\subset M$. For $V\in\{V_u,V_v\}$ the change-of-variables formula for this inverse branch gives
\[
 \vol_M(j(V\setminus C))=\int_{V\setminus C}|\operatorname{Jac}j|\,d\vol_U=0,
\]
since $\xi(U\setminus C)=0$ and $\vol_U\ll\xi$. Choose
\[
        y\in W\setminus\bigl(j(V_u\setminus C)\cup
                              j(V_v\setminus C)\bigr).
\]
Then $(j|_{V_u})^{-1}(y),(j|_{V_v})^{-1}(y)\in C$ are distinct since $V_u\cap V_v=\varnothing$. This contradicts injectivity of $j|_C$. Hence $\#j^{-1}(y)=1$ for $y\in j(U)$. The local inverses therefore define $j^{-1}\in C^\infty(j(U),U)$.
\end{proof}

Let $(P,\eta)$ be a second countable analytic $d$-manifold with finitely many components and a smooth positive probability density. Let $\mathfrak l$ be a real Lie algebra of dimension $d$, and let $\beta_p:T_pP\to\mathfrak l$ be an analytic coframe satisfying $[X_u,X_v]=X_{[u,v]}$, where $X_u(p)=\beta_p^{-1}u$ for $u\in\mathfrak l$. Let $\mathcal G$ be a second countable Lie group with Lie algebra $\mathfrak l$. Assume it acts on the right on the measure algebra of $(P,\eta)$ (Definition~\ref{H-def:probability-algebra-action}). Thus $R_{gh}=R_h R_g$ and, for every measurable $A\subset P$,
\[
 \eta(R_g[A])=\eta(A),\qquad \lim_{g\to e}\eta(R_g[A]\mathbin\triangle[A])=0 .
\]

Choose a basis $u_1,\ldots,u_d$ of $\mathfrak l$. For $w:\mathbb R^d\to\mathcal G$, $w(t)=\exp(t_1u_1)\cdots\exp(t_du_d)$, $D_0w(e_i)=u_i$. Thus $w$ is an analytic chart on a sufficiently small ball about $0$. If $\phi_u^s$ is the local flow of $X_u$, define $\Phi(p,t)=\phi_{u_d}^{t_d}\circ\cdots\circ\phi_{u_1}^{t_1}(p)$ on the open set of $(p,t)\in P\times\mathbb R^d$ for which every partial composition exists. It takes values in $P$ and $D_t\Phi(p,0)e_i=X_{u_i}(p)$.

\begin{lemma}[Realization from local flow words]\label{H-lem:coframe-measured-realization}
Let a countable group $\Gamma$ preserve $\eta$ ergodically. Suppose $\sigma:\Gamma\to\operatorname{Aut}(\mathcal G)$ is a homomorphism with $\gamma R_g\gamma^{-1}=R_{\sigma_\gamma(g)}$ on the measure algebra for $\gamma\in\Gamma$, $g\in\mathcal G$. Suppose that on each local word chart the measured action of $w(t)$ agrees with $\Phi(p,t)$ for $(\eta\otimes\mathcal L^d)$-almost every $(p,t)$. In a point realization (a Borel action inducing the prescribed action on the measure algebra), the union of finitely many $\mathcal G$-orbits has full measure. The normalizing group $\Gamma$ permutes them transitively. Fix one such orbit $O$, let $\Gamma_O=\{\gamma\in\Gamma:\gamma O=O\}$ and put $s=[\Gamma:\Gamma_O]$. Then $O$ is isomorphic as a right $\mathcal G$-space to $\Delta\backslash\mathcal G$ for a lattice $\Delta<\mathcal G$. There is an analytic open embedding
\[
 F:P\longrightarrow Y_{\rm fr}=\{1,\ldots,s\}\times O
\]
agreeing $\eta$-almost everywhere with the point realization. Its image is conull for the measure assigning mass $1/[\Gamma:\Gamma_O]$ to each copy of $O$, with Haar probability as normalized restriction. For $p\in P$ and $u\in\mathfrak l$, $F(\phi_u^t(p))=F(p)\cdot\exp(tu)$ for all $t$ in the existence interval of the local flow at $p$.
\end{lemma}
\begin{proof}
\textbf{The measured orbits.} As in~\cite{KRHArithmetic}, compare the local flow coordinates with the orbits of a point action. Use the left action of $\Gamma$ on $P$ and the right action of $\mathcal G$. Their normalization relation is $\gamma(p\cdot g)=(\gamma p)\cdot\sigma_\gamma(g)$ modulo $\eta$, for fixed $\gamma\in\Gamma$ and $g\in\mathcal G$. In the right action of $S=\mathcal G\rtimes_\sigma\Gamma$, the element $(e,\gamma)$ therefore acts as $p\mapsto\gamma^{-1}p$. Give $\Gamma$ the discrete topology. This is a second countable locally compact group. On the standard probability algebra each open coset acts continuously in the group parameter by strong continuity of the $\mathcal G$-action. By Mackey's point-realization theorem~\cite[Theorem~1, p.~330]{Mackey} there are a standard Borel right $S$-space $\Omega$ and a measure-space isomorphism modulo null sets $\mathcal F:(P,\eta)\to(\Omega,\mathcal F_*\eta)$. For each fixed $g\in S$, it intertwines $R_g$ with right translation at $\eta$-almost every $p\in P$. Choose a Borel conull $P_0$ on which $\mathcal F$ is injective. The local-word hypothesis and Fubini give, for $\eta$-almost every $p\in P$,
\[
 \Phi(p,t)\in P_0,\qquad \mathcal F(\Phi(p,t))=\mathcal F(p)\cdot w(t)
 \quad(\mathcal L^d\text{-a.e. }t\text{ in the local word domain}).
\]
Here $\Phi(p,\cdot)$ is a local diffeomorphism near $0$, so its parameter law there is absolutely continuous with respect to $\eta$. Fix such a $p$ and take an identity neighborhood $U$ on which $k\mapsto\Phi(p,w^{-1}(k))$ is injective. For $d\in U\cap\Delta_p$, where $\Delta_p=\{g\in\mathcal G:\mathcal F(p)\cdot g=\mathcal F(p)\}$, the word endpoints for $k$ and $dk$ lie in $P_0$ and have equal $\mathcal F$-images for Haar-almost every $k\in U\cap d^{-1}U$. This open intersection contains the identity. Injectivity first of $\mathcal F|_{P_0}$ and then of the word gives $dk=k$, hence $d=e$. Thus $\Delta_p$ is discrete and closed.
The stabilizers satisfy $\operatorname{Stab}_{\mathcal G}(\mathcal F(p)\cdot g)=g^{-1}\Delta_pg$ and $\operatorname{Stab}_{\mathcal G}(\gamma\mathcal F(p))=\sigma_\gamma(\Delta_p)$ for a normalizer with automorphism $\sigma_\gamma$.
The orbit $\{\mathcal F(p)\cdot g:g\in\mathcal G\}$ is Borel because its parameter map from
$\Delta_p\backslash \mathcal G$ is an injective Borel map between standard
Borel spaces~\cite[Corollary~15.2]{Kechris}. For each $p$ in the preceding $\eta$-full set the local coframe word covers an open patch $U_p\subset P$, with $\eta(U_p)>0$ and $\eta\{y\in U_p:\mathcal F(y)\notin\mathcal F(p)\mathcal G\}=0$. Thus $(\mathcal F_*\eta)(\mathcal F(p)\mathcal G)\ge\eta(U_p)>0$. The restriction of $\mathcal F_*\eta$ to this orbit is a nonzero finite right-invariant measure. Pull it back by $\Delta_p g\mapsto\mathcal F(p)\cdot g$. Uniqueness of invariant quotient measure gives a scalar multiple of Haar quotient measure on $\Delta_p\backslash\mathcal G$, with the scalar fixed by the mass of the orbit~\cite{WeilIntegration}. For $p\in P$ as above, let $w(t)=\exp(t_1u_1)\cdots\exp(t_du_d)$ as before the lemma. Its local endpoint satisfies
\begin{equation}\label{H-eq:realization-word}
                   \mathcal F(\Phi(p,t))=\mathcal F(p)\cdot w(t)
\end{equation}
for $\mathcal L^d$-almost every $t$ in a neighborhood of $0\in\mathbb R^d$. Both $t\mapsto\Phi(p,t)$ and
$t\mapsto\Delta_p w(t)$ are now local analytic diffeomorphisms.
They identify a neighborhood of $p$ analytically with an open subset of $\Delta_p\backslash\mathcal G$. For $\eta$-almost every $y$ in that neighborhood the point of $\Omega$ represented by this chart is $\mathcal F(y)$.

\par\Needspace{3\baselineskip}\smallskip\noindent\textbf{An analytic representative at every point.} For any $x\in P$, the local coframe word $\Phi(c,t)$ gives a local diffeomorphism $(c,t)\mapsto(c,\Phi(c,t))$ near $(x,0)$. Its image contains a product about $(x,x)$. Shrink so the group word is a coordinate chart throughout. Choose $c$ in the first factor of that product so that \eqref{H-eq:realization-word} holds for $\mathcal L^d$-almost every $t$. Such $c$ have full $\eta$-measure in the factor by Fubini. Solve $\Phi(c,t)=x$. Since $\Delta_c$ is discrete, the model word is a local diffeomorphism also at this $t$. This gives an analytic representative near $x$.
Take a countable subcover $(U_i)_{i\in I}$, $I\subset\mathbb N$, of these charts. Write $F_i:U_i\to O_i$ for the chart into its orbit. Under the orbit parameterization, $\eta\{p\in U_i:F_i(p)\ne\mathcal F(p)\}=0$. If $U_i\cap U_j\ne\varnothing$, then $\eta(U_i\cap U_j)>0$. Hence some $p$ there satisfies $\mathcal F(p)\in O_i\cap O_j$, which implies $O_i=O_j$. Fix $o$ in this common orbit and set $\Delta_o=\{g\in\mathcal G:o\cdot g=o\}$. When $\mathcal F(c)=o\cdot a$, the chart in $\Delta_o\backslash\mathcal G$ is $\Phi(c,t)\mapsto\Delta_o a w(t)$. Equation~\eqref{H-eq:realization-word} identifies both $F_i(p)$ and $F_j(p)$ with $\mathcal F(p)$ for $\eta$-almost every $p\in U_i\cap U_j$. The disagreement set is open by continuity and has zero measure. Full support of $\eta$ makes it empty. The charts glue to a local diffeomorphism $F$. The orbit containing $F(p)$ is locally constant in $p\in P$, so each connected component meets just one of the orbits $O_i$.

\par\Needspace{3\baselineskip}\smallskip\noindent\textbf{The finite family and its measure.} Let $\mathscr O=\{O_i:i\in I\}$ be the resulting family of orbits in $\Omega$. The preceding argument gives $|\mathscr O|\le|\pi_0(P)|<\infty$. Every $O'\in\mathscr O$ has positive mass. The normalizers permute these orbits by $\gamma(z\mathcal G)=(\gamma z)\mathcal G$. Ergodicity makes $\Gamma$ transitive on $\mathscr O$, since the union of the translates of any one orbit is invariant and has positive mass. Fix an orbit $O\subset\Omega$ in this family and $o\in O$. Writing $\Delta=\{g\in\mathcal G:o\cdot g=o\}$ gives $O=\Delta\backslash\mathcal G$. The orbit measure computed above gives finite~covolume. Write $\Gamma=\bigsqcup_{j=1}^{s}t_j\Gamma_O$. Then $(j,y)\mapsto t_jy$ identifies $\{1,\ldots,s\}\times O$ with the union of the frame orbits. In these coordinates the right action is $(j,y)\cdot g=(j,y\cdot\sigma_{t_j}^{-1}(g))$, since $(t_jy)\cdot g=t_j(y\cdot\sigma_{t_j}^{-1}(g))$. Each copy has mass $1/s$. Write $m_{\rm fr}$ for this probability on $Y_{\rm fr}$.

Choose a Borel set $C\subset P_0$ of full $\eta$-measure on which $F=\mathcal F$. Then $F|_C$ is injective. Apply Lemma~\ref{H-lem:conull-local-diffeomorphism} to $F:P\to Y_{\rm fr}$ with this $C\subset P$ and $\xi=\eta\sim\vol_P$. Its local analytic inverses give an analytic open embedding. Its measured representative has law $m_{\rm fr}$, so $F_*\eta=m_{\rm fr}$ and $m_{\rm fr}(Y_{\rm fr}\setminus F(P))=0$. On each open word domain the maps $(p,t)\mapsto F(\Phi(p,t))$ and $(p,t)\mapsto F(p)\cdot w(t)$ are continuous. By the local-word hypothesis they agree for $\eta\otimes\mathcal L^d$-almost every $(p,t)$ there. This measure has full support, so equality holds throughout the domain. Substitute $t=se_i$ and differentiate at $s=0$. Linearity in $u$ gives $D_pF\,X_u(p)=\left.\frac{d}{dt}\right|_0 F(p)\exp(tu)$ for all $u\in\mathfrak l$. Uniqueness of flows gives the identity for $\phi_u^t$.
\end{proof}

\begin{proof}[Proof of Theorem~\ref{main}]
\leavevmode\hypertarget{proof-main-start}{}Apply Lemma~\ref{H-lem:coframe-measured-realization} with $P=P_1$, $\beta=\vartheta$ the analytic coframe of Lemma~\ref{coframe}, $\mathcal G=L$ the group of Lemma~\ref{adjoin}, and the countable group of normalizers generated by $\Lambda$ and the $T_t$ for a finite set of representatives of $\Gamma_1/\Lambda$. These normalizers preserve $\eta$, and $\Lambda$ acts ergodically. For two normalizers $T,S$ and $g\in L$,
\[
 R_{\sigma_{TS}(g)}=(TS)R_g(TS)^{-1}
     =T R_{\sigma_S(g)}T^{-1}=R_{\sigma_T(\sigma_S(g))}.
\]
The $L$-action is faithful by Lemma~\ref{adjoin}, so $\sigma_{TS}=\sigma_T\circ\sigma_S$, as required. For a fixed tuple of times, Theorem~\ref{H-thm:coframe-integration} identifies each measured factor with its local flow modulo $\eta$. Each preceding factor preserves $\eta$, so its inverse image of the exceptional set is null. Fubini gives agreement of the measured word and the local flow word for $\eta\otimes\mathcal L^{\dim P_1}$-almost every $(p,t)$ in its domain, as required by the lemma.

\par\Needspace{3\baselineskip}\smallskip\noindent\textbf{The quotient by $K$.} Write $F_{\rm fr}:P_1\to Y_{\rm fr}$ for the resulting open frame embedding. Fix a normalizer $T$ with induced automorphism $\alpha$ of $L$. On a model orbit with basepoint $o$, put $\Delta=\{g\in L:og=o\}$. Choose a basepoint $o'$ in the orbit containing $T(o)$ and put $\Delta'=\{g\in L:o'g=o'\}$. Choose $b\in L$ with $T(o)=o'b$. Then the induced map between these orbits is
\[
 \Delta g\longmapsto\Delta' b\alpha(g),\qquad b\alpha(\Delta)b^{-1}=\Delta',
\]
because $T(og)=T(o)\alpha(g)$ for the point action. This map is analytic, so the almost-everywhere equivariance of $F_{\rm fr}$ holds everywhere by continuity and full support. For each fixed $h\in K$ the same argument gives $F_{\rm fr}(ph)=F_{\rm fr}(p)h$ everywhere. Let $q:Y_{\rm fr}\to Y_1=Y_{\rm fr}/K$ be the quotient map on the frame model of Lemma~\ref{H-lem:coframe-measured-realization}. Since the right $K$-action on $P_1$ is free and $F_{\rm fr}$ is injective, $K$ acts freely on the open image $F_{\rm fr}(P_1)$. Set
\[
 Y_1^\circ=q(F_{\rm fr}(P_1)),\qquad
 i_1:Y_1^\circ\longrightarrow O_1',\quad i_1(q(F_{\rm fr}(p)))=c(p).
\]
For a local analytic section $\sigma:U\to P_1$ over $U\subset O_1'$ from Lemma~\ref{sections}, the equivariant analytic diffeomorphism
\[
 U\times K\longrightarrow F_{\rm fr}(c^{-1}(U)),\qquad
 (x,h)\longmapsto F_{\rm fr}(\sigma(x)h)=F_{\rm fr}(\sigma(x))h
\]
identifies the quotient projection with $(x,h)\mapsto x$. Thus $Y_1^\circ$ lies in the smooth locus, and $i_1^{-1}(x)=q(F_{\rm fr}(\sigma(x)))$ on $U$. Two sections differ by right multiplication in $K$ and give the same inverse. The normalizers preserve $K$ and their composition defects are right $K$-translations by Lemma~\ref{finite}. After quotienting, the group law therefore holds on the open image. That image has full measure and is dense, so continuity extends the law to all of $Y_1$. The open sets $Y_1^\circ$ and $O_1'$ are invariant under $\Gamma_1$. With $\nu_1=q_* (F_{\rm fr})_*\eta$ we have $i_1q F_{\rm fr}=c$ on $P_1$, so
\[
 (i_1)_*(\nu_1|_{Y_1^\circ})=m_1,\qquad m_1|_{O_1'}=k\,m|_{O_1'} .
\]

\par\Needspace{3\baselineskip}\smallskip\noindent\textbf{Returning to the original lattice.} Recall from \eqref{components} that $m=k^{-1}\sum_{j=1}^{k}m_j$, with $m_j=k\,m|_{X_j}$, where the $m_j$ are the distinct $\Gamma$-translates of $m_1$. Here $O_1'=\bigcap_{t\in\Gamma_1}tO_1$ is the open set of Lemma~\ref{finite}, with $\Gamma_1=\{t\in\Gamma:t_*m_1=m_1\}$. Choose representatives $r_1=e,\ldots,r_k$ of $\Gamma/\Gamma_1$. Distinct $(r_j)_*m_1$ are mutually singular. If $r_jO_1'\cap r_lO_1'$ were nonempty for $j\ne l$, this open set would have positive $m$-measure and both component measures would restrict there to $k\,m$, a contradiction. Put
\[
 \begin{gathered}
 Y=\Gamma\times_{\Gamma_1}Y_1,\qquad Y^\circ=\Gamma\times_{\Gamma_1}Y_1^\circ,\qquad
 \nu=\frac1k\sum_{j=1}^k(y\mapsto[r_j,y])_*\nu_1,\\
 O=\bigsqcup_{j=1}^k r_jO_1',\qquad
 i:Y^\circ\longrightarrow O,\quad i([r_j,y])=r_j i_1(y)\quad(y\in Y_1^\circ).
 \end{gathered}
\]
On each copy $i$ is an analytic diffeomorphism onto $r_jO'_1$. For $\gamma\in\Gamma$ write $\gamma r_j=r_lh$ with $h\in\Gamma_1$. Then $i(\gamma[r_j,y])=r_li_1(hy)=\gamma r_ji_1(y)$. Since $\nu_1(Y_1^\circ)=1$,
\[
 \nu(Y^\circ)=1,\qquad
 i_*(\nu|_{Y^\circ})
 =\frac1k\sum_{j=1}^k(r_j)_*(i_1)_*(\nu_1|_{Y_1^\circ})
 =\frac1k\sum_{j=1}^k(r_j)_*m_1=m .
\]
Thus $m(O)=1$ and $\overline O=M$ by full support. Since $i:Y^\circ\to O$ intertwines the actions and sends $\nu|_{Y^\circ}$ to the ergodic probability $m$, $\nu$ is ergodic.

Each connected component of $Y$ is open. Indeed the images of the cosets of $L^0$ in a double quotient are connected and open, and two such images are either equal or disjoint. The model probability has full support, so ergodicity makes the union of the translates of a component conull. Its complement is open of zero mass and hence empty. All components have the same positive mass, hence their number $s$ is finite. Choose $C$ in the copy $[e,Y_1]$ and put $\Gamma_C=\{\gamma\in\Gamma:\gamma C=C\}$. Since $\Gamma_C\le\Gamma_1$, the number of connected components is $s=[\Gamma:\Gamma_C]=k[\Gamma_1:\Gamma_C]$ and $[\gamma,y]\mapsto\gamma y$ identifies $\Gamma\times_{\Gamma_C}C$ with $Y$.

\par\Needspace{3\baselineskip}\smallskip\noindent\textbf{The connected model.} To present $C$ as a double quotient, choose the basepoint in $C$, conjugating its lattice if necessary. In its model summand $\Delta\backslash L/K$, put $L'=\Delta L^0$ and $K'=L'\cap K$. The image of $L^0$ under $L\to\Delta\backslash L/K$ is exactly $C$, by the component description above. Since $\Delta L^0=L'$, the map
\[
 \Delta\backslash L'/K'\longrightarrow C,\qquad \Delta gK'\longmapsto\Delta gK
\]
is onto. If $g,g'\in L'$ and $g'=\delta g h$ with $\delta\in\Delta$, $h\in K$, then $h=g^{-1}\delta^{-1}g'\in L'\cap K=K'$, which gives injectivity. The local quotient charts agree since $L'$ is open in $L$. The open subset $\Delta\backslash L'\subset\Delta\backslash L$ has finite positive Haar mass, so $\Delta$ is a lattice in $L'$.

Fix $\gamma\in\Gamma_C$. In the summand containing $C$, its affine representative $g\mapsto b\alpha(g)$ satisfies $\alpha(K)=K$ and $b\alpha(\Delta)b^{-1}=\Delta$. Since $\Delta bK\in C$, write $b=\delta l h$ with $\delta\in\Delta$, $l\in L^0$, $h\in K$. Put $\alpha'=\operatorname{Ad}(h)\circ\alpha$. Then
\[
 \Delta b\alpha(g)K=\Delta l\alpha'(g)K,\qquad
 \alpha'(\Delta)=l^{-1}\Delta l,\qquad
 \alpha'(L')=l^{-1}\Delta lL^0=\Delta L^0=L'.
\]
Also $\alpha'(K)=K$, hence $\alpha'(K')=K'$. Thus $g\mapsto l\alpha'(g)$ induces an affine map on $\Delta\backslash L'/K'$ in the sense of Section~\ref{sec:homogeneous}. To identify the probability on $C$, let $\lambda$ be Haar probability on $\Delta\backslash L$ in this summand, put $A=\Delta\backslash L'$, and write $q(\Delta g)=\Delta gK$. Since $K'=L'\cap K$ is open in the compact group $K$, the coset space $K'\backslash K$ is finite. We have
\[
 q^{-1}(C)=\bigsqcup_{K'h\in K'\backslash K}A h,\qquad q(zh)=q(z)\quad(z\in A,\ h\in K).
\]
Right invariance of $\lambda$ gives $(q_*\lambda)(B)=[K:K']\,\lambda(A\cap q^{-1}(B))$ for every Borel $B\subset C$. Substituting $B=C$ and dividing by $(q_*\lambda)(C)=[K:K']\,\lambda(A)$ shows that the normalized restriction to $C$ is the image of $\lambda|_A/\lambda(A)$, Haar probability on $\Delta\backslash L'$. Taking $L'$ as the model group and $H=K'$, the measure on the $s=[\Gamma:\Gamma_C]$ copies of $C$ is therefore $s^{-1}\sum_{j=1}^s\delta_j\otimes m_C$, as asserted.
\end{proof}

\paragraph{Acknowledgements.}
We thank Amir Algom, Snir Ben Ovadia, Aaron Brown, Alex Eskin and Simion Filip for helpful discussions. Federico Rodriguez Hertz acknowledges support from NSF grant DMS-2453688 and from the Anatole Katok Chair in Mathematics. Zhiren Wang was partially supported by NSF grant DMS-2453689.
\paragraph{Use of artificial intelligence.}
The authors devised the proof strategy and its division into cases. At the same time, recent AI tools have provided the technical force that allowed to check all the branching cases. AI assisted with the proofs of technical lemmas and with the writing. The authors take full responsibility for this manuscript.

\end{document}